\documentclass[10pt, a4paper]{amsart}
\usepackage{amsfonts}
\usepackage{amsthm}
\usepackage{amsmath}
\usepackage{longtable}
\usepackage{amssymb}
\usepackage{latexsym}
\usepackage{hyperref}

\theoremstyle{plain} \newtheorem{lem}{Lemma}[section]
\theoremstyle{plain} \newtheorem{prop}[lem]{Proposition}
\theoremstyle{plain} 
\theoremstyle{plain} 

\theoremstyle{definition}

\newtheorem{rem}[lem]{Remark}

\numberwithin{equation}{section}

\newcommand{\C}{\mathbb{C}}

\newcommand{\R}{\mathbb{R}}
\newcommand{\HH}{\mathbb{H}}

\newcommand{\GL}{\mathrm{\mathop{GL}}}
\newcommand{\End}{\mathrm{\mathop{End}}}

\newcommand{\ad}{\mathrm{\mathop{ad}}}

\newcommand{\Aut}{\mathrm{\mathop{Aut}}}

\newcommand{\ssl}{\mathfrak{\mathop{sl}}}
\newcommand{\su}{\mathfrak{\mathop{su}}}
\renewcommand{\u}{\mathfrak{\mathop{u}}}
\newcommand{\so}{\mathfrak{\mathop{so}}}
\renewcommand{\sp}{\mathfrak{\mathop{sp}}}
\newcommand{\g}{\mathfrak{g}}
\newcommand{\h}{\mathfrak{h}}
\renewcommand{\d}{\mathfrak{d}}

\newcommand{\uu}{\mathfrak{u}}

\newcommand{\kk}{\mathfrak{k}}
\newcommand{\pp}{\mathfrak{p}}
\newcommand{\mya}{\mathfrak{a}}

\newcommand{\A}{\mathcal{A}}
\newcommand{\veps}{\varepsilon}

\begin{document}

\title[Maximal reductive subalgebras]
{Real forms of embeddings of maximal reductive subalgebras of the complex
simple Lie algebras of rank up to 8}

\author{Willem A. de Graaf}
\address{Dipartimento di Matematica\\
Universit\`{a} di Trento\\
Italy}
\email{degraaf@science.unitn.it}

\author{Alessio Marrani}
\address{Museo Storico della Fisica e Centro Studi e Ricerche Enrico Fermi,
  Roma, Italy, and Dipartimento di Fisica e Astronomia Galileo Galilei,
  Universit\`a di Padova, and INFN, sezione di Padova, Italy}
\email{alessio.marrani@pd.infn.it}

\begin{abstract}
  We give tables of noncompact
  real forms of maximal reductive subalgebras of complex
  simple Lie algebras of rank up to 8. These were obtained by computational
  methods that we briefly describe. We also discuss applications in theoretical
  physics of these embeddings.
\end{abstract}

\maketitle

\section{Introduction}

The classification of maximal reductive subalgebras of
exceptional and classical complex Lie algebras
is well known. A systematic study started with the work of Dynkin
\cite{dyn2, dyn}. Subsequently, it was developed in a number of works, such
as \cite{LG, Minchenko, DG1}. For what
concerns non-compact real forms, the maximal reductive
subalgebras of exceptional and classical Lie algebras which give rise to
symmetric embeddings are listed e.g. in Table 9.7 of \cite{Gilmore},
which reports results taken from a rather vast literature (see the references
therein). On the other hand, the non-compact real forms of maximal,
semisimple (reductive) subalgebras of exceptional Lie algebras which give
rise to non-symmetric embeddings are studied in a number of papers scattered
in works of mathematics and of theoretical physics, which are not easily
collected; for instance, in the high-energy physics literature
(Maxwell-Einstein-scalar theories, possibly endowed with some amount of
local supersymmetry), some maximal subalgebras giving rise to non-symmetric
embeddings have been considered in \cite{ADFMT-1, super-Ehlers, FMZ-D=5}.
In the mathematical literature Gray \cite{Gray} classified the real forms of
the maximal reductive subalgebras that arise as fixed point subalgebras of
automorhisms of order three (overlooking two cases, see Section
\ref{sec:ehlers}).
Komrakov \cite{Kom} also classified all real forms of the
maximal (non-symmetric) embedding $A_{4}\oplus A_{4}$ $\subset E_{8}$, and he
also carried out a classification of the maximal S-subalgebras of
exceptional real Lie algebras.

Recently algorithms, together with their implementations, have been developed
to compute the real forms of embeddings of complex semisimple Lie algebras,
\cite{dfg2, fg}. However, a classificatory and exhaustive approach
was not pursued in these papers. In the present paper, by developing and
systematically applying these techniques we aim at obtaining the complete
list of noncompact
real forms of embeddings of maximal reductive subalgebras of the
simple complex Lie algebras of ranks up to 8.

Some remarks to clarify our aim
are in order. Firstly, with a real form of an embedding of complex Lie algebras
we mean a pair of real forms $\mathfrak{a},\g$ of complex reductive Lie algebras
$\mathfrak{a}^c, \g^c$ such that $\mathfrak{a}^c\subset \g^c$ and $\mathfrak{a}
\subset \g$. Secondly, it is possible to have a pair of complex semisimple
(or reductive) Lie algebras $\mathfrak{a}^c\subset \g^c$ with
real forms $\mathfrak{a},\g$ such that $\mathfrak{a}$ is maximal in $\g$
without $\mathfrak{a}^c$ being maximal in $\g^c$. In fact, Komrakov \cite{Kom}
has a list of cases where that happens. We do not quite understand
this list (for example, the cases appearing under item a) in the cited paper
seem to be
maximal reductive subalgebras, and also appear in our tables), nor do we
have methods to detect those cases independently. Therefore we restrict to
listing the real forms of embeddings of the maximal complex reductive
subalgebras of complex semisimple Lie algebras.

The second part of the paper
has the tables that we obtained containing the maximal reductive subalgebras
of the noncompact simple real Lie algebras of rank up to 8. These subalgebras
are also contained, in explicit form, in the latest version of the
{\sf CoReLG} package for {\sf GAP}4 (\cite{corelg}).

The remainder of this paper is divided in a few sections. In Section
\ref{sec:maxred} we describe the maximal reductive subalgebras of the
simple complex Lie algebras. Section \ref{sec:phys} gives an overview of
some applications of real embeddings of reductive Lie algebras in high-energy
theoretical physics.
After that we have Section \ref{sec:methods} giving a brief overview of the
computational methods that we used to obtain real forms of complex embeddings.
Finally Section \ref{sec:tables} has the tables of those real forms.

\section{Maximal reductive subalgebras of complex semisimple Lie algebras}
\label{sec:maxred}

In this section we describe the complex embeddings of which we compute the
real forms.

Let $\g^c$ be a complex simple Lie algebra.
After Dynkin \cite{dyn}, a subalgebra $\mya^c\subset \g^c$ is said to
be {\em regular} if there is a Cartan subalgebra $\h^c$ of $\g^c$ with
$[\h^c,\mya^c]\subset \mya^c$. Such a subalgebra is spanned by root spaces
with respect to $\h^c$ along with $\h^c\cap \mya^c$. A subalgebra
of $\g^c$ is called an {\em R-subalgebra} if it is contained
in a proper regular subalgebra of $\g^c$. A subalgebra that is not an
R-subalgebra is called an S-subalgebra.

Dynkin showed that every non-semisimple subalgebra of $\g^c$ is an R-subalgebra
(\cite{dyn}, Theorem 7.3). Hence S-subalgebras are necessarily
semisimple. So the maximal reductive S-subalgebras coincide with the
maximal semisimple S-subalgebras. These have been classified by Dynkin
\cite{dyn2,dyn}, and are also contained in the lists obtained in \cite{DG1}.

Now we consider the regular maximal reductive subalgebras. Let
$\Delta$ be a set of simple roots of the root system of $\g^c$ with respect
to a Cartan subalgebra $\h^c\subset \g^c$. Furthermore, let $\delta$
denote the highest root of the root system and write
$$\delta = \sum_{\alpha\in \Delta} n_\alpha \alpha$$
where the $n_\alpha$ are positive integers.

Let $e_\alpha$, $f_\alpha$ (for $\alpha$ in the root system of $\g^c$)
be elements spanning the root spaces corresponding
to $\alpha$ and $-\alpha$ respectively. Let $h_\alpha$ be a scalar multiple
of $[e_\alpha,f_\alpha]$ such that $[h_\alpha,e_\alpha] = 2e_\alpha$.
Let $\alpha\in \Delta$. Then by $\g^c(\alpha)$ we denote the subalgebra
generated by $e_\beta,f_\beta$ for $\beta\in \Delta\setminus \{\alpha\}$ along
with $e_\delta$, $f_\delta$. By $\g^c[\alpha]$ we denote the subalgebra
generated by $e_\beta,f_\beta$ for $\beta\in \Delta\setminus \{\alpha\}$ along
with $e_\alpha$, $h_\alpha$. Then $\g^c(\alpha)$ is a maximal subalgebra if
and only if $n_\alpha$ is prime (this statement goes back to \cite{borelsieb},
see also \cite{gotogros}, Chapter 8). 
The $\g^c(\alpha)$ are semisimple, so they are maximal reductive
subalgebras. On the other hand, a maximal reductive, non semisimple subalgebra
is contained in a $\g^c[\alpha]$. The latter has a unique (up to conjugation)
maximal reductive
subalgebra, which we denote by $\g^c[\alpha]'$; it is generated by
$e_\beta,f_\beta$ for $\beta\in \Delta\setminus \{\alpha\}$ along with $h_\alpha$.
However, it is possible that $\g^c[\alpha]'$ is contained in a $\g^c(\beta)$.
This situation can be characterized as follows. Consider the Dynkin diagram of
the set of roots $\Delta\cup\{-\delta\}$ (this is called the extended Dynkin
diagram). Now $\g^c[\alpha]'$ is contained in a $\g^c(\beta)$ if and only if
the diagram obtained by removing the node corresponding to $\alpha$ is not
the Dynkin diagram of $\g^c$. (Indeed, if the latter diagram is not the
Dynkin diagram of $\g^c$ then $\g^c(\alpha)$ is not equal to $\g^c$ and
$\g^c[\alpha]'\subset \g^c(\alpha)$. Conversely, if the mentioned diagram is
the Dynkin diagram of $\g^c$ then none of the $\g(\beta)$ contain the semisimple
part of $\g^c[\alpha]'$.) By inspection this then leads to the maximal
reductive R-subalgebras listed in Table \ref{tab:maxred}.

\begin{table}[htb]

  \begin{tabular}{|c|c|c|c|}
    \hline
    type & max reductive & type & max reductive\\
    \hline
    $A_n$ & $A_k + A_{n-1-k}+T_1$ & $D_n$ & $A_{n-1}+T_1$\\
    & $1\leq k\leq n-1$ &        & $D_{n-1}+T_1$\\
    \hline
    $B_n$ & $B_{n-1}+T_1$ & $E_6$ & $D_5+T_1$ \\
    \hline
    $C_n$ & $A_{n-1}+T_1$ & $E_7$ & $E_6+T_1$\\
    \hline
  \end{tabular}
  \caption{Maximal reductive non-semisimple subalgebras of the simple Lie
    algebras. Here $T_1$ denotes a $1$-dimensional centre.}\label{tab:maxred}
\end{table}

\begin{rem}
  The subalgebras of Table \ref{tab:maxred} are also given by Dynkin in
  \cite{dyn}, Table 12a. This table lists the regular semisimple subalgebras
  $\mya^c\subset \g^c$ such that no semisimple regular subalgebra $\hat\mya^c$
  exists with $\mya^c \subsetneq \hat\mya^c \subsetneq \g^c$. However, in one
  case this appears to be not quite correct. In fact, for the embeddings
  $B_{n-1}\subset B_n$ we have the chain
  $$B_{n-1} \subset D_n \subset B_n.$$
  Here both subalgebras $B_{n-1}$ and $D_n$ are regular, but they are not
  normalized by the same Cartan subalgebra.

  We also have the dual chain $D_{n-1}\subset B_{n-1} \subset D_n$. But here,
  obviously, $B_{n-1}$ is not regular in $D_n$.
\end{rem}

\begin{rem}
  Here we point out two oversights that are present in some places in the
  literature. Firstly, in various works (such as
  \cite{McKay-Patera, Slansky, Yamatsu}), the subalgebra
  $A_{1}\oplus A_{1}\oplus A_{1}$ is reported to be maximal in $D_{6}$, while
  actually it is not. We have that $A_{1}\oplus A_{1}\oplus A_{1}$ is a maximal
  (non-symmetric) subalgebra in $A_{1}\oplus C_{3}$, which in turn is maximal
  (and non-symmetric) in $D_{6}$. This has been noted in the tables in
  \cite{LG}.

  Secondly, the maximal S-subalgebra $C_{3}$
in $C_{7}$ (giving rise to a non-symmetric embedding) is not listed in some
works (such as \cite{McKay-Patera, Slansky, Yamatsu}),
whereas it is considered in \cite{dyn2}, \cite{LG}, Table VII.
All in all, the existence of a maximal $C_{3}$ in $C_{7}$ is a consequence
of the anti-self-conjugation (i.e., symplecticity) of the irreducible
representation $\mathbf{14}^{\prime }=\wedge _{0}^{3}\mathbf{6}$ of $C_{3}$.
It is here worth
recalling that the action of $C_{3}$ on the $\mathbf{14}^{\prime }$ is
\textquotedblleft of type $E_{7}$" \cite{brown}. In \cite{Kac} it has
been proved that such action has a finite number of nilpotent orbits, with
one-dimensional ring of invariant polynomials generated by a
quartic homogeneous polynomial. The latter is related to the square of the
Bekenstein-Hawking entropy of extremal black hole solutions
to the \textquotedblleft magic" Maxwell-Einstein $\mathcal{N}=2$
supergravity having the split real form $\mathfrak{sp}(6,\mathbb{R})$ of
$C_{3}$ as electric-magnetic duality symmetry \cite{GST} (see \cite{AM-Rev}
for a review and a list of references).
\end{rem}

\section{Physical applications of real embeddings of reductive Lie algebras}
\label{sec:phys}

\subsection{Super-Ehlers embeddings, and their non-supersymmetric versions}
\label{sec:ehlers}

The non-compact real form $\mathfrak{sl}(5,\mathbb{R})\oplus
\mathfrak{sl}(5,\mathbb{R})\subset E_{8(8)}$ of the maximal
(non-symmetric)
embedding $A_{4}\oplus A_{4}$ $\subset E_{8}$ is known in supergravity as an
example of \textquotedblleft super-Ehlers" embedding, concerning the
maximally supersymmetric Einstein gravity in 7 space-time dimensions.
Super-Ehlers embeddings, which unify the Ehlers gravity embeddings with the
global electric-magnetic duality symmetries of Einstein-Maxwell theories (at
least in the cases with symmetric scalar manifolds), have been introduced
and studied, in presence of underlying (local) supersymmetry, in \cite%
{super-Ehlers}; in particular, in the Appendices of \cite{super-Ehlers} a
general proof of existence of such regular and rank-preserving embeddings,
which are non-symmetric in most cases, is given, within an approach
completely different from the ones employed in the present paper. The
general structure of super-Ehlers embeddings is the following (with $%
3\leqslant D\leqslant 11$ denoting the number of Lorentzian space-time
dimensions):%
\begin{equation}
\mathfrak{g}_{D,\mathcal{N}}\oplus \mathfrak{sl}(D-2,\mathbb{R})\subset
\mathfrak{g}_{3,\mathcal{N}};
\end{equation}%
where $\mathfrak{g}_{D,\mathcal{N}}$ is the electric-magnetic duality Lie algebra of $D$-dimensional Maxwell-Einstein
theories endowed with $2\mathcal{N}$ local supersymmetries (corresponding,
in $D=3$, to $\mathcal{N}$-extended supergravity), and $\mathfrak{sl}(D-2,%
\mathbb{R})$ is the Ehlers symmetry Lie algebra in $D$ Lorentzian
dimensions. In presence of supersymmetry, super-Ehlers embeddings are listed
and classified in \cite{super-Ehlers} (for the $D=5$ case, see also \cite%
{FMZ-D=5}). Note that
\begin{eqnarray*}
E_{6(6)}\oplus \mathfrak{sl}(3,\mathbb{R}) &\subset &E_{8(8)},\\
E_{6(-26)}\oplus \mathfrak{sl}(3,\mathbb{R}) &\subset &E_{8(-24)}
\end{eqnarray*}
are super-Ehlers
embeddings for $\mathcal{N}=16$ and $\mathcal{N}=4$ supergravity theories in
$D=5$, respectively. (These two embeddings have been overlooked in
\cite{Gray}.)
On the other hand, non-supersymmetric Maxwell-Einstein
theories (coupled to non-linear sigma model of scalar fields) in various
dimensions are not considered in \cite{super-Ehlers}, but nevertheless they
display some \textquotedblleft non-supersymmetric Ehlers embeddings"; whose
some examples list as follows :%
\begin{eqnarray}
\mathfrak{so}(6,6)\oplus \mathfrak{sl}(2,\mathbb{R}) &\subset &E_{7(7)}; \\
\mathfrak{sl}(6,\mathbb{R})\oplus \mathfrak{sl}(3,\mathbb{R}) &\subset
&E_{7(7)}; \\
\mathfrak{sl}(6,\mathbb{R})\oplus \mathfrak{sl}(2,\mathbb{R}) &\subset
&E_{6(6)}; \\
\mathfrak{sl}(3,\mathbb{R})\oplus \mathfrak{sl}(3,\mathbb{R})\oplus
\mathfrak{sl}(3,\mathbb{R}) &\subset &E_{6(6)}.
\end{eqnarray}%
These embeddings have been discussed in the so-called non-supersymmetric
magic Maxwell-Einstein theories in \cite{MPRR} (see also \cite{Romano}), in
which the role of the Ehlers symmetry and related truncations has been
highlighted.

\subsection{Cubic Jordan algebras, their symmetries and related embeddings}

Various embeddings reported in the tables of the present paper have an
interpretation in terms of symmetries of cubic Jordan algebras. Such
symmetries are given by the $\mathfrak{der}$ (derivations), $\mathfrak{str}%
_{0}$ (reduced structure), $\mathfrak{conf}$ (conformal) and $\mathfrak{qconf%
}$ (quasi-conformal) Lie algebras associated to a given cubic Jordan algebra. The quasiconformal realizations of
non-compact groups were discovered by the authors of \cite{Gunaydin:2000xr} and were further developed
and applied in \cite{Gunaydin:2005mx,Gun-2,Gunaydin:2009zza,Gunaydin:2007qq,Gunaydin:2007bg}. The conformal groups associated with Jordan algebras
were studied much earlier in \cite{Gunaydin:1975mp,Gunaydin:1992zh}. Over the reals $\mathbb{R}$, the Lie algebras $\mathfrak{str}_{0}$%
, $\mathfrak{conf}$ and $\mathfrak{qconf}$ respectively correspond to the
electric-magnetic duality ($U$-duality\footnote{%
Here $U$-duality is referred to as the \textquotedblleft continuous"
symmetries of \cite{j1, j2}. Their discrete versions are the $U$-duality
non-perturbative string theory symmetries introduced in \cite{HT}.}) Lie
algebras of some Maxwell-Einstein supergravity theories in $D=5,4,3$
Lorentzian space-time dimensions (cfr. e.g. \cite{Gun-2}, and the references
therein;
in $D=3$ all vectors need to be dualized into scalars); such symmetries are
non-linearly realized on the scalars, while vectors do sit in some linear
representations of them. Jordan algebras, such as $J_{3}^{\mathbb{H}%
_{s}}$, $J_{3}^{\mathbb{C}_{s}}$ and $\mathbb{R}\oplus \mathbf{\Gamma} _{m,n}$ with $m
$ (or $n$) $\neq 1$ and $5$, are associated to non-supersymmetric models
\cite{MPRR, Romano}. When considering simple cubic Jordan algebras, all the
aforementioned related Lie algebras fill the first ($\mathfrak{der}$),
second ($\mathfrak{str}_{0}$), third ($\mathfrak{conf}$) and fourth ($%
\mathfrak{qconf}$) rows of the relevant Magic Square of
Freudenthal-Rozenfeld-Tits \cite{FRT1, FRT2, FRT3}, and the following maximal
embeddings hold (see
\cite{Squaring-Magic}, and references therein):%
\begin{eqnarray}
\mathfrak{der} &\subset &\mathfrak{str}_{0};  \label{s1} \\
\mathfrak{der}\oplus \mathfrak{sl}(2,\mathbb{R}) &\subset &\mathfrak{conf};
\label{s1-bis} \\
\mathfrak{str}_{0}\oplus \mathfrak{so}(1,1) &\subset &\mathfrak{conf};
\label{s2} \\
\mathfrak{conf}\oplus \mathfrak{sl}(2,\mathbb{R}) &\subset &\mathfrak{qconf};
\label{s3} \\
\mathfrak{str}_{0}\oplus \mathfrak{sl}(3,\mathbb{R}) &\subset &\mathfrak{%
qconf}; \\
\mathfrak{der}\oplus G_{2(2)} &\subset &\mathfrak{qconf}.  \label{s3-tris}
\end{eqnarray}%
\begin{table}[h]
\begin{center}
\begin{tabular}{|c|c|c|c|c|}
\hline
& $\mathbb{R}$ & $\mathbb{C}$ & $\mathbb{H}$ & $\mathbb{O}$ \\ \hline
$\mathbb{R}$ & $SO(3)$ & $SU(3)$ & $USp(3)$ & $F_{4(-52)}$ \\ \hline
$\mathbb{C}_s$ & $SL(3,\mathbb{R})$ & $SL(3,\mathbb{C})$ & $SL(3,\mathbb{H})$ & $%
E_{6(-26)}$ \\ \hline
$\mathbb{H}_s$ & $Sp(6,\mathbb{R})$ & $SU(3,3)$ & $SO^*(12)$ & $E_{7(-25)}$
\\ \hline
$\mathbb{O}_s$ & $F_{4(4)}$ & $E_{6(2)}$ & $E_{7(-5)}$ & $E_{8(-24)}$ \\
\hline
\end{tabular}%
\end{center}
\caption{The \textit{single-split} MS ${\mathcal{L}}_{3}(\mathbb{A}_{s},%
\mathbb{B})$ \protect\cite{GST} (see \protect\cite{Squaring-Magic} for
details)}
\label{tab:MS1}
\end{table}
\begin{table}[h]
\begin{center}
\begin{tabular}{|c|c|c|c|c|}
\hline
& $\mathbb{R}$ & $\mathbb{C}_s$ & $\mathbb{H}_s$ & $\mathbb{O}_s$ \\ \hline
$\mathbb{R}$ & $SO(3)$ & $SL(3,\mathbb{R})$ & $Sp(3,\mathbb{R})$ & $F_{4(4)}$
\\ \hline
$\mathbb{C}_s$ & $SL(3,\mathbb{R})$ & $SL(3,\mathbb{R}) \times SL(3,\mathbb{R%
})$ & $SL(6,\mathbb{R})$ & $E_{6(6)}$ \\ \hline
$\mathbb{H}_s$ & $Sp(6,\mathbb{R})$ & $SL(6,\mathbb{R})$ & $SO(6,6)$ & $%
E_{7(7)}$ \\ \hline
$\mathbb{O}_s$ & $F_{4(4)}$ & $E_{6(6)}$ & $E_{7(7)}$ & $E_{8(8)}$ \\ \hline
\end{tabular}%
\end{center}
\caption{The \textit{double-split} MS ${\mathcal{L}}_{3}(\mathbb{A}_{s},%
\mathbb{B}_{s})$ \protect\cite{BS} (see \protect\cite{Squaring-Magic} for
further details)}
\label{tab:MS2}
\end{table}
Within the physical interpretation of such symmetries as $U$-duality Lie
algebras, the commuting $\mathfrak{so}(1,1)$ in (\ref{s2}) can be regarded
as the Kaluza-Klein compactification radius of the $S^{1}$-reduction $%
D=5\rightarrow 4$; alternatively, such an $\mathfrak{so}(1,1)$ can also be
conceived as the Lie algebra of the pseudo-K\"{a}hler connection of the
pseudo-special K\"{a}hler (and pseudo-Riemannian) symmetric coset\footnote{%
As intuitive, the names starting with a capital letter denote the Lie groups
whose Lie algebra is the same name starting lowercase. Moreover,
\textquotedblleft $mcs$" denotes the maximal compact subgroup.} $\frac{Conf}{%
Str_{0}\otimes SO(1,1)}$, obtained from $\frac{Str_{0}}{mcs(Str_{0})}$ by
applying the inverse $R^{\ast }$-map pertaining to a timelike
compactification $D=s+t=4+1\rightarrow 4+0$ where $s$ and $t$ respectively
denote the number of spacelike and timelike dimensions
\cite{R-map-1, R-map-2, timelike-reduction}.
On the other hand, the commuting $\mathfrak{sl}(2,\mathbb{R})$ in (\ref{s3})
can be identified with the Ehlers symmetry arising from the $S^{1}$%
-reduction $D=4\rightarrow 3$; such an $\mathfrak{sl}(2,\mathbb{R})$ can
also be interpreted as the connection of the para-quaternionic (and
pseudo-Riemannian) symmetric coset $\frac{QConf}{Conf\otimes SL(2,\mathbb{R})%
}$, obtained from $\frac{Conf}{mcs\left( Conf\right) }$ by applying the
inverse $c^{\ast }$-map pertaining to a timelike compactification $%
D=(3,1)\rightarrow (3,0)$ \cite{BGM, c-map, timelike-reduction}. As the
embeddings (\ref{s1})-(\ref{s3-tris}) are obtained by moving along the
columns of the relevant (real form of the) Magic Square (for a fixed row
entry), another class of embeddings can be obtained by moving along the rows
of the relevant Magic Square (for a fixed column entry). In the symmetric
(real forms of the) rank-3 Magic Square, as the double-split ${\mathcal{L}}%
_{3}(\mathbb{A}_{s},\mathbb{B}_{s})$ given in Table \ref{tab:MS2}, these
embeddings trivially coincide with (\ref{s1})-(\ref{s3-tris}), but their
intepretation corresponds to the restriction from one (division $\mathbb{A}$
or split $\mathbb{A}_{s}$) algebra to a smaller one. For the non-symmetric,
single-split Magic Square ${\mathcal{L}}_{3}(\mathbb{A}_{s},\mathbb{B})$
reported in Table \ref{tab:MS1}, different maximal
embeddings hold true, which are
different from the ones given in (\ref{s1})-(\ref{s3-tris}); namely
($G_2^{\mathrm{cpt}}\equiv G_{2(-14)}$) :%
\begin{eqnarray}
  \d(J_{3}^{\mathbb{R}}) &\subset \d(J_{3}^{\mathbb{C}}) \\
  \su(2)\oplus \d(J_{3}^{\mathbb{R}}) &\subset \d(J_3^{\mathbb{H}})\\
  \mathfrak{u}(1)\oplus \d(J_{3}^{\mathbb{C}}) &\subset \d(J_3^{\mathbb{H}})\\
  \su(2)\oplus \d(J_{3}^{\mathbb{H}}) &\subset \d(J_3^{\mathbb{O}})\\
  \su(3)\oplus \d(J_{3}^{\mathbb{C}}) &\subset \d(J_3^{\mathbb{O}})\\
  G_2^{\mathrm{cpt}}\oplus \d(J_{3}^{\mathbb{R}}) &\subset \d(J_3^{\mathbb{H}})
\end{eqnarray}
(where $\d$ can be each of $\mathfrak{der}$, $\mathfrak{str}_0$,
$\mathfrak{conf}$, $\mathfrak{qconf}$).

Similar embeddings holds for simple Lorentzian cubic Jordan algebras (see
\cite{Squaring-Magic}).

\subsection{Semisimple subalgebras of simple Jordan algebras, and their symmetries}

Another remarkable class of embedding stems from the relation
between simple cubic Jordan algebras \cite{JVNW} $J_{3}^{\mathbb{A}}$ or $J_{3}^{\mathbb{%
A}_{s}}$ and some elements of the (bi-parametric) infinite sequence of
semi-simple Jordan algebras $\mathbb{R}\oplus \mathbf{\Gamma }_{m,n}$
mentioned above, exploiting the Jordan-algebraic isomorphisms $J_{2}^{%
\mathbb{A}}\cong \mathbf{\Gamma }_{1,q+1}\left( \cong \mathbf{\Gamma }%
_{q+1,1}\right) $ and $J_{2}^{\mathbb{A}_{s}}\cong \mathbf{\Gamma }%
_{q+2+1,q/2+1}$, where $q:=$dim$_{R}\mathbb{A}=8,4,2,1$ for $\mathbb{A}=%
\mathbb{O},\mathbb{H},\mathbb{C},\mathbb{R}$, and $q:=$dim$_{R}\mathbb{A}%
_{s}=8,4,2$ for $\mathbb{A}_{s}=\mathbb{O}_{s},\mathbb{H}_{s},\mathbb{C}_{s}$
(see \textit{e.g.} App. A of \cite{Gun-2} - and Refs. therein - for an
introduction to division and split algebras). Indeed, the following
(maximal, rank-preserving) Jordan-algebraic embeddings hold :%
\begin{eqnarray}
J_{3}^{\mathbb{A}} &\supset &\mathbb{R}\oplus J_{2}^{\mathbb{A}}\cong
\mathbb{R}\oplus \mathbf{\Gamma }_{1,q+1}; \\
J_{3}^{\mathbb{A}_{s}} &\supset &\mathbb{R}\oplus J_{2}^{\mathbb{A}%
_{s}}\cong \mathbb{R}\oplus \mathbf{\Gamma }_{q/2+1,q/2+1}.
\end{eqnarray}%
Thus, one can consider their consequences at the level of symmetries of
cubic Jordan algebras defined over the corresponding algebras, obtaining :%
\begin{eqnarray}
  \d(\mathbb{R}\oplus J_{2}^{\mathbb{A}})\oplus \mathcal{A}_{q} &\subset
  \d(J_{3}^{\mathbb{A}}) \label{embs!-1} \\
  \d(\mathbb{R}\oplus J_{2}^{\mathbb{A}_{s}}) \oplus
    \widetilde{\mathcal{A}}_{q} &\subset \d(J_{3}^{\mathbb{A}_{s}})\label{embs!-8}
\end{eqnarray}
(where as above $\d$ can be each of $\mathfrak{der}$, $\mathfrak{str}_0$,
$\mathfrak{conf}$, $\mathfrak{qconf}$).

Note the maximal nature of the embeddings (\ref{embs!-1})-(\ref{embs!-8}),
as well as the presence of the commuting algebras $\mathcal{A}_{q}$ and $%
\widetilde{\mathcal{A}}_{q}$, defined as follows:%
\begin{eqnarray}
\mathcal{A}_{q} &:&=\mathfrak{tri}\left( q\right) \ominus so(q)=\varnothing ,%
\mathfrak{so}(3),\mathfrak{so}(2),\varnothing \text{ for~}q=8,4,2,1; \\
\widetilde{\mathcal{A}}_{q} &:&=\widetilde{\mathfrak{tri}}\left( q\right)
\ominus \widetilde{so}(q)=\varnothing ,\mathfrak{sl}(2,\mathbb{R}),\mathfrak{%
so}(1,1)\text{ for~}q=8,4,2,
\end{eqnarray}%
where $\mathfrak{tri}\mathbb{\ }$and $\mathfrak{so}$ respectively denote the
triality and orthogonal (norm-preserving) symmetries (and they are tilded
when pertaining to split algebras; see \textit{e.g. }\cite{CFMZ1-D=5, CCM},
and Refs. therein). Analogous results hold for Lorentzian Jordan algebras.
Within the physical ($U$-duality) interpretation, $\mathcal{A}_{q}$ and $%
\widetilde{\mathcal{A}}_{q}$ are consistent with the properties of spinors
in $q+2$ dimensions, with Lorentzian signature $(t,s)=\left( 1,q+1\right) $
respectively
Kleinian signature $(t,s)=\left( q/2+1,q/2+1\right) $; indeed, the
electric-magnetic ($U$-duality) symmetry Lie algebra in $D=6$ (Lorentzian)
space-time dimensions is $\mathfrak{so}\left( 1,q+1\right) \oplus \mathcal{A}%
_{q}$ for $\mathbb{A}$-based theories (which are endowed with minimal,
chiral $\left( 1,0\right) $ supersymmetry) and $\mathfrak{so}\left(
q/2+1,q/2+1\right) \oplus \widetilde{\mathcal{A}}_{q}$ for $\mathbb{A}_{s}$%
-based theories (which are non-supersymmetric for $q=2,4$, and endowed with
maximal, non-chiral $\left( 2,2\right) $ supersymmetry for $q=8$); \textit{%
cfr. e.g.} \cite{Kugo-Townsend}, \cite{magic-D=6} (and Refs. therein) and
\cite{MNY} for further discussion.

\section{Computational methods}\label{sec:methods}

In this section we briefly describe the computational methods that we used
to construct the real forms of the maximal semisimple (or reductive)
subalgebras of the simple Lie algebras of ranks up to 8.
We have one procedure for constructing the regular subalgebras (taken from
\cite{dfg2}) and one procedure for the S-subalgebras (from \cite{fg}).
First we recall some general notation, and subsequently describe the
methods in two subsections. Our main reference for the general theory
is \cite{onishchik}.

By $\iota\in \C$ we denote the imaginary unit.
Let $\g^c$ be a complex simple Lie algebra. An {\em anti-involution} of $\g^c$
is a map $\eta : \g^c\to \g^c$ with $\eta(x+y)=\eta(x)+\eta(y)$, $\eta(\mu x)
=\bar \mu x$, $\eta([x,y])=[\eta(x),\eta(y)]$, $\eta(\eta(x))=x$ for all
$x,y\in \g^c$, $\mu\in \C$. If $\eta : \g^c\to \g^c$ is any map then we
set $\g^c_\eta = \{ x\in \g^c \mid \eta(x)=x\}$.

A real form $\g$ of $\g^c$ is given by three maps, $\tau,\sigma,\theta : \g^c
\to \g^c$ where
\begin{itemize}
\item $\tau,\sigma$ are anti-involutions and $\theta$ is an involution,
\item $\theta=\tau\sigma=\sigma\tau$,
\item $\g^c_\tau$ is compact and $\g^c_\sigma = \g$,
\item setting $\kk = \{ x\in \g \mid \theta(x)=x \}$,
  $\pp = \{ x\in \g \mid \theta(x)=-x\}$; then $\g = \kk
  \oplus \pp$ is a Cartan decomposition of $\g$ (note that $\theta
  (\g)=\g$ because $\theta$ commutes with $\sigma$); the restriction of
  $\theta$ to $\g$ is called a {\em Cartan involution} of $\g$,
\item set $\u=\g^c_\tau$, then $\theta$ leaves $\u$ invariant so that
  $\u = \u_1\oplus \u_{-1}$ (the eigenspaces of $\theta$ with eigenvalues
  $\pm 1$) and $\kk = \u_1$, $\pp = \iota \u_{-1}$.
\end{itemize}

Let $\mya\subset \g$ be a semisimple subalgebra and $\mya^c =
\mya+\iota \mya$, which is a semisimple subalgebra of $\g^c$.
Then $\mya$ is a real form of the complex subalgebra $\mya^c$.
So $\mya$ has a Cartan decomposition $\mya =
\kk_{\mya} \oplus \pp_{\mya}$. It follows from the Karpelevich-Mostow
theorem (\cite{onishchik}, Corollary 1 of \S 6) that $\mya$ is conjugate
by an inner automorphism to a subalgebra $\mya'$ such that
$\kk_{\mya'} \subset \kk$ and $\pp_{\mya'} \subset \pp$. Equivalently, $\theta$
maps $\mya'$ to itself, and its restriction to $\mya'$ is a Cartan
involution of $\mya'$. So we may restrict to finding $\theta$-stable subalgebras
of $\g$.

\subsection{Real forms of regular subalgebras}\label{alg:reg}

We say that a subalgebra $\mya$ of $\g$ is
regular if it is normalized by a Cartan subalgebra $\h$ of $\g$. If we want
to stress the particular Cartan subalgebra that we are referring to we
say that $\mya$ is $\h$-regular. In that case
$\mya^c$ is spanned by root spaces of $\g^c$ (relative to $\h^c$) along with
$\mya^c\cap \h^c$. By $\Psi(\h^c,\mya^c)$ we denote the set of roots involved
in this. It is a subset of the root system of $\g^c$ with respect to $\h^c$.

Let $G^c$ denote the adjoint group of $\g^c$. This group can be characterized
in several ways. Firstly it is the connected algebraic subgroup of $\GL(\g^c)$
with Lie algebra $\ad \g^c$. Secondly it is the group generated by
all elements $\exp( ad x )$ where $x\in \g^c$ is nilpotent. Thirdly, it is
the identity component of $\Aut(\g^c)$. By $G$ we denote the group consisting
of all $g\in G^c$ such that $g(\g)=\g$. If we represent elements of $G^c$
by their matrices with respect to a basis of $\g$ then $G=G^c(\R)$, the
group of all elements of $G^c$ with coefficients in $\R$. Let $G_0$ be the
identity component of $G$. In this section we outline how to obtain the
regular semisimple subalgebras of $\g$ up to conjugacy by $G_0$.

First of all, Dynkin (\cite{dyn}) devised an algorithm to obtain the regular
semisimple subalgebras of $\g^c$ up to conjugacy by $G^c$ (see \cite{gra16},
\S 5.9 for a recent account). We can use this algorithm to obtain all
$\h^c$-regular semisimple subalgebras up to $G^c$-conjugacy, for a fixed
Cartan subalgebra $\h^c$ of $\g^c$. We have that two $\h^c$-regular semisimple
subalgebras $\mya_1^c$, $\mya_2^c$ of $\g^c$ are $G^c$-conjugate if and only
if $\Psi(\h^c,\mya_1^c)$, $\Psi(\h^c,\mya_2^c)$ are $W(\g^c,\h^c)$-conjugate,
where the latter denotes the Weyl group of the root system of $\g^c$
with respect to $\h^c$.

If one is only interested in subalgebras of $\g^c$ then it suffices to
consider just one Cartan subalgebra as they are all conjugate under
$G^c$. In general the real form $\g$ has more Cartan subalgebras that
are non-conjugate. Sugiura (\cite{sugiura}) proved that $\g$ has a finite
number of Cartan subalgebras up to $G_0$-conjugacy. In \cite{dfg}
Sugiura's method was made into an algorithm for listing the Cartan
subalgebras of $\g$ up to conjugacy by $G_0$.

Let $\mya$ be a regular semisimple subalgebra of $\g$. Then the normalizer
$\mathfrak{n}_\g(\mya) = \{ x\in \g \mid [x,\mya]\subset \mya\}$ is reductive.
Therefore it has a {\em unique} maximally noncompact Cartan subalgebra
(that is, a Cartan subalgebra whose intersection with $\pp$ has maximal
dimension, see \cite{knapp02}, Proposition 6.61). We say that $\mya$ is
{\em strongly $\h$-regular} if $\h$ is a maximally noncompact Cartan
subalgebra of $\mathfrak{n}_\g(\mya)$.

Let $\h$ be a Cartan subalgebra of $\g$. Set $N = \{g\in G_0 \mid
g(\h)=\h\}$ and $Z=\{ g\in G_0 \mid g(x)=x \text{ for all } x\in \h\}$.
Then $W(\g,\h) = N/Z$ is called the {\em real Weyl group} of $\h$.
It can naturally be identified with a subgroup of $W(\g^c,\h^c)$, (see
\cite{knapp02}, \S VII.8). There are algorithms to compute $W(\g,\h)$,
see \cite{realweyl}.

We have the following criterion: Let $\mya_1,\mya_2$ be two strongly
$\h$-regular semisimple subalgebras of $\g$. They are $G_0$-conjugate if and
only if $\Psi(\h^c,\mya_1^c)$, $\Psi(\h^c,\mya_2^c)$ are conjugate under
$W(\g,\h)$ (see \cite{dfg2}, Proposition 24). This gives an immediate method
for deciding whether $\mya_1$, $\mya_2$ are $G_0$-conjugate.

Now the algorithm for obtaining all regular semisimple subalgebras of $\g$
(up to $G_0$-conjugacy) runs as follows:

\begin{enumerate}
\item Compute the Cartan subalgebras of $\g$ (up to $G_0$-conjugacy). For each
  obtained Cartan subalgebra $\h$ perform the following steps:
  \begin{enumerate}
  \item Compute the $\h^c$-regular subalgebras of $\g^c$ using Dynkin's
    algorithm. For each obtained subalgebra $\mya^c$ perform the following
    steps:
    \begin{enumerate}
    \item Compute the stabilizer $S$ of $\Psi(\h^c,\mya^c)$ in $W=W(\g^c,
      \h^c)$. Compute a list of representatives $w_1,\ldots,w_s$ of the
      double cosets $W(\g,\h)w_i S$ in $W$.
    \item Construct the $\h^c$-regular subalgebras $\mya_i^c$ of $\g^c$
      with $\Psi(\h^c,\mya_i^c) = w_i\cdot \Psi(\h^c,\mya^c)$.
    \item Throw away the $\mya_i^c$ that are not $\sigma$-stable.
    \item Of the remaining ones compute a basis of $\mya_i=\mya_i^c\cap \g$, and
      throw away those that are not strongly $\h$-regular. Add the remaining
      ones to the final list.
    \end{enumerate}
  \end{enumerate}
\end{enumerate}

As before, let $\theta$ denote the Cartan involution of $\g$.
The Cartan subalgebras found in the first step are $\theta$-stable (see the
algorithm in \cite{dfg}).
Therefore the subalgebras constructed by the above procedure are automatically
$\theta$-stable as well (\cite{dfg2}, Proposition 21).

Note that the above algorithm is correct. Indeed, each $\h^c$-regular semisimple
subalgebra $\mathfrak{b}^c$ of $\g^c$ that is $G^c$-conjugate to $\mya^c$
has $\Psi(\h^c,\mathfrak{b}^c) =
w\cdot \Psi(\h^c,\mya^c)$ for some $\mya^c$ from the initial list.
Furthermore, note that the $w_i\cdot \Psi(\h^c,\mya^c)$ exhaust the images
under $W(\g^c,\h^c)$ of $\Psi(\h^c,\mya^c)$ up to conjugacy by
$W(\g,\h)$.

\subsection{Real forms of S-subalgebras}\label{alg:S}

From Section \ref{sec:maxred} we recall that a subalgebra of $\g^c$ is called an
{\em S-subalgebra} if it is not contained in a proper regular subalgebra.
For these subalgebras we have no method to list all real
forms in $\g$ up to conjugacy by $G_0$. Therefore we are more modest and
consider the following question. Let $\veps : \g^c\hookrightarrow \tilde\g^c$
be an embedding of semisimple complex Lie algebras and let $\g$ be a real form
of $\g^c$.
{\em Find (up to isomorphism) the real forms $\tilde\g$ of $\tilde\g^c$ such
  that $\veps( \g ) \subset \tilde\g$.}
Because two real forms $\tilde\g,\tilde\g'$ of $\tilde\g^c$ are isomorphic if
and only if there is a $\phi\in \Aut(\tilde\g^c)$ with $\phi(\tilde\g)=
\tilde\g'$ (\cite{onishchik}, \S 2,
Proposition 1) we may replace the given embedding
$\veps$ by $\phi\circ \veps$ for any $\phi\in \Aut(\tilde\g^c)$.

Let $\theta,\tilde\theta$ be Cartan involutions of $\g$, $\tilde\g$
respectively. A classical theorem of Karpelevich states that $\veps(\g)
\subset \tilde \g$ if and only if $\veps\circ \theta = \tilde \theta \circ
\veps$ (see \cite{onivin}, Theorem 3.7 of Chapter 6). Here we consider
$\theta$ and $\g$ to be given, and hence have to construct real forms
$\tilde \g$ with Cartan involution $\tilde\theta$ satisfying the mentioned
condition. Our method is based on Proposition \ref{prop:eps}.

Let $\sigma : \g^c\to \g^c$ be the conjugation with respect to the real form
$\g$. Let $\u$ be a fixed compact form of $\g^c$ with corresponding conjugation
$\tau : \g^c\to \g^c$ with $\sigma\tau=\tau\sigma$. Set $\theta = \sigma\tau$.
Then the restriction of $\theta$ to $\g$ is a Cartan involution of $\g$.
From \cite{fg} we have the following result.

\begin{prop}\label{prop:eps}
Let $\tilde\g\subset \tilde\g^c$ be a real form of $\tilde\g^c$ such that
$\veps(\g)\subset \tilde \g$. Then there are a compact form $\tilde\uu\subset \tilde\g^c$
of $\tilde\g^c$, with conjugation $\tilde\tau : \tilde\g^c\to\tilde\g^c$, and an involution
$\tilde\theta$ of $\tilde\g^c$ such that
\begin{enumerate}
\item \label{eps1}  $\veps(\uu)\subset \tilde \uu$,
\item \label{eps2} $\veps\theta = \tilde\theta\veps$,
\item \label{eps3} $\tilde\theta \tilde\tau = \tilde\tau\tilde\theta$,
\item \label{eps4} there is a Cartan decomposition $\tilde\g = \tilde\kk \oplus \tilde \pp$,
such that the restriction of $\tilde\theta$ to $\tilde\g$
is the corresponding Cartan involution, and
$\tilde\uu = \tilde \kk\oplus \imath \tilde\pp$.
\end{enumerate}
Conversely, if $\tilde\uu\subset \tilde\g$ is a compact form, with corresponding conjugation
$\tilde\tau$, and $\tilde\theta$ is an involution of $\tilde\g^c$ such that \eqref{eps1},
\eqref{eps2} and \eqref{eps3} hold,
then $\tilde\theta$ leaves $\tilde\uu$ invariant, and setting $\tilde\kk = \tilde\uu_1$,
$\tilde\pp = \imath\tilde\uu_{-1}$ (where $\tilde\uu_k$ is the $k$-eigenspace of $\tilde\theta$),
we get that $\tilde\g = \tilde\kk\oplus\tilde\pp$ is a real form of $\tilde\g^c$ with
$\veps(\g)\subset \tilde\g$.
\end{prop}

To construct the real forms $\tilde\g$ the strategy is now the following.
First we fix a compact form $\tilde\u$ of $\tilde\g^c$ and replace $\veps$ by
$\phi\veps$ for a $\phi\in \Aut(\tilde\g^c)$ in order to have $\veps(\u)\subset
\tilde\u$. Secondly we construct the space
$$\A = \{ A\in \End(\tilde\g^c) \mid A (\ad(\veps\theta(y))) = (\ad\veps(y))A
\text{ for all } y\in \g^c\}.$$
(Note that a basis of $\A$ can be computed by solving a set of linear
equations; \cite{fg} also has a more efficient method to construct $\A$.)
Let $\tilde\theta$ be an involution of $\tilde\g^c$. Then $\veps\theta =
\tilde\theta\veps$ if and only if $\tilde\theta\in \A$ and
\begin{equation}\label{eq:theta}
\tilde\theta
(\ad x) \tilde\theta= \ad \tilde\theta(x) \text{ for all } x\in \tilde\g^c
\end{equation}
(\cite{fg}, Proposition 3.6).
Hence $\A$ contains the maps $\tilde\theta$. The
conditions that $\tilde\theta$ be an involution, \eqref{eq:theta}, \eqref{eps3}
are translated to polynomial equations on the coefficients of $\tilde\theta$
with respect to a basis of $\A$. These polynomial equations are then
studied and solved by means of the technique of Gr\"obner bases
(see \cite{clo}).

For the details of this procedure we refer to \cite{fg}. But we do remark that
if $\veps(\g)\subset\tilde \g$ for some real form $\tilde \g$ of $\tilde\g^c$
then we get just one subalgebra of $\tilde\g$. However, there may be more
subalgebras $\mya \subset \tilde \g$ such that $\mya^c$ is $G^c$-conjugate to
$\veps(\g^c)$, without $\mya$ being $G_0$-conjugate to $\veps(\g)$. This method
does not detect such a situation (unlike the method for the regular case).
Furthermore, there are cases in type $D_n$ for $n$ even, where there are
subalgebras $\mya_1^c,\mya_2^c\subset \tilde\g^c$ that are not conjugate under
$G^c$, but are conjugate under an outer automorphism. In these cases we
just consider one of the corresponding embeddings in $\tilde \g^c$ because
changing the embedding $\veps$ to a $\phi\veps$ may identify the subalgebras
$\mya_1^c$, $\mya_2^c$.

\subsection{Implementation}

We have implemented these methods in the computer algebra system {\sf GAP}4
(\cite{gap4}) using the packages {\sf SLA} (\cite{sla}) and {\sf CoReLG}
(\cite{corelg}). The implementation of the algorithm of Section
\ref{alg:reg} is quite straightforward. For the procedure indicated in
Section \ref{alg:S} we remark that the package {\sf SLA} contains
tables with the semisimple subalgebras of the simple complex Lie algebras.
The embeddings $\veps$ are simply given by those tables.

\section{The tables}\label{sec:tables}

Here we give the tables that we have obtained. In these we use the
notation of real forms as used in \cite{knapp02}. Several authors use
different notation, for instance, $\mathfrak{su}^{\ast }(2n)$ instead of
$\mathfrak{sl}(4,\mathbb{H})$. Moreover, we use the mathematical
notation for symplectic algebras : $C_{n}$ corresponds to $\sp_{n}$, and not
to $\sp_{2n}$, as in a large part of the physics literature.
We also denote by $G_2^\mathrm{cpt}$, $F_4^\mathrm{cpt}$, $E_6^\mathrm{cpt}$ and
$E_7^\mathrm{cpt}$ the compact real forms of $G_2$, $F_4$, $E_6$ and $E_7$,
respectively. In the literature, they are also indicated by $G_{2(-14)}$,
$F_{4(-52)}$, $E_{6(-78)}$ and $E_{7(-133)}$.

As written in the previous section, the algorithm for the regular subalgebras
determines them up to conjugacy by the adjoint group. Therefore, occasionally
it happens that we get more than one real embedding of the same type.
This is indicated in the tables by putting the number of those embeddings
(e.g., as $2\times$). For the S-subalgebras we can only determine whether
there is an embedding, except when the S-subalgebra is of type $A_1$. In those
cases the subalgebras correspond to the nilpotent orbits in $\g$
(see \cite{colmcgov}); so for those subalgebras we also know how many there
are up to conjugacy by the adjoint group. It may also happen that there is more
than one S-subalgebra of type $A_1$. In those cases the different occurrences
of these subalgebras are distinguished by an upper index, giving the
Dynkin index of the subalgebra.

The S-subalgebras of the real forms of the exceptional algebras have also
been determined by Komrakov \cite{Kom}. However, his table for $E_6$ has
quite a lot of omissions (for example, the real forms $E_{6(6)}$, $E_{6(-26)}$
appear to have been forgotten altogether). Also the table for the S-subalgebras
of the real forms of type $E_6$ in \cite{fg} have a number of errors and
omissions (some subalgebras have been identified incorrectly, others have been
omitted).

\vspace{3mm}
\noindent{\bf Type $A_2$}\\

\begin{longtable}{|l|c|c|}
\caption{Maximal subalgebras of $\su(1,2)$.}\label{tab:su12}
\endfirsthead
\hline
\multicolumn{3}{|l|}{\small\slshape subalgebras of $\su(1,2)$} \\
\hline
\endhead
\hline
\endfoot
\endlastfoot

\hline
type & real subalgebras & \\
\hline

$A_1$ & $
   \begin{aligned}
\ssl(2,\R)\oplus\u(1)\\
\su(2)\oplus\u(1)\\
   \end{aligned} $ & R \\
   \hline

$A_1$ & $\ssl(2,\R)$ & S\\
   \hline

\end{longtable}

\begin{longtable}{|l|c|c|}
\caption{Maximal subalgebras of $\ssl(3,\R)$.}\label{tab:sl3R}
\endfirsthead
\hline
\multicolumn{3}{|l|}{\small\slshape subalgebras of $\ssl(3,\R)$} \\
\hline
\endhead
\hline
\endfoot
\endlastfoot

\hline
type & real subalgebras & \\
\hline

$A_1$ & $\ssl(2,\R)\oplus\so(1,1)$ & R\\
   \hline

$A_1$ & $
   \begin{aligned}
\su(2)\\
\ssl(2,\R)\\
   \end{aligned} $ & S \\
   \hline

\end{longtable}

\noindent{\bf Type $B_2$}

\begin{longtable}{|l|c|c|}
\caption{Maximal subalgebras of $\so(2,3)$.}\label{tab:so23}
\endfirsthead
\hline
\multicolumn{3}{|l|}{\small\slshape subalgebras of $\so(2,3)$} \\
\hline
\endhead
\hline
\endfoot
\endlastfoot

\hline
type & real subalgebras & \\
\hline

$A_1+A_1$ & $
   \begin{aligned}
\ssl(2,\R)\oplus \ssl(2,\R)\\
\ssl(2,\C)\\
   \end{aligned} $ & R \\
   \hline

$A_1+T_1$ & $
   \begin{aligned}
\ssl(2,\R)\oplus\so(1,1)\\
\ssl(2,\R)\oplus\u(1)\\
\su(2)\oplus\u(1)\\
   \end{aligned} $ & R \\
   \hline

$A_1$ & $\ssl(2,\R)$ & S\\
   \hline

\end{longtable}

\begin{longtable}{|l|c|c|}
\caption{Maximal subalgebras of $\so(4,1)$.}\label{tab:so41}
\endfirsthead
\hline
\multicolumn{3}{|l|}{\small\slshape subalgebras of $\so(4,1)$} \\
\hline
\endhead
\hline
\endfoot
\endlastfoot

\hline
type & real subalgebras & \\
\hline

$A_1+A_1$ & $
   \begin{aligned}
\ssl(2,\C)\\
\su(2)\oplus \su(2)\\
   \end{aligned} $ & R \\
   \hline

   $A_1+T_1$ & $
   \begin{aligned}
\ssl(2,\R)\oplus\u(1)\\
\su(2)\oplus\so(1,1)\\
   \end{aligned} $ & R \\
   \hline

\end{longtable}

\noindent{\bf Type $G_2$}

\begin{longtable}{|l|c|c|}
\caption{Maximal subalgebras of $G_{2(2)}$.}\label{tab:G22}
\endfirsthead
\hline
\multicolumn{3}{|l|}{\small\slshape subalgebras of $G_{2(2)}$} \\
\hline
\endhead
\hline
\endfoot
\endlastfoot

\hline
type & real subalgebras & \\
\hline

$A_2$ & $
   \begin{aligned}
\ssl(3,\R)\\
\su(1,2)\\
   \end{aligned} $ & R \\
   \hline

$A_1+A_1$ & $
   \begin{aligned}
\ssl(2,\R)\oplus \ssl(2,\R)\\
\su(2)\oplus \su(2)\\
   \end{aligned} $ & R \\
   \hline

$A_1$ & $\ssl(2,\R)$ & S\\
   \hline

\end{longtable}

\noindent{\bf Type $A_3$}

\begin{longtable}{|l|c|c|}
\caption{Maximal subalgebras of $\su(1,3)$.}\label{tab:su13}
\endfirsthead
\hline
\multicolumn{3}{|l|}{\small\slshape subalgebras of $\su(1,3)$} \\
\hline
\endhead
\hline
\endfoot
\endlastfoot

\hline
type & real subalgebras & \\
\hline

$A_2$ & $
   \begin{aligned}
\su(1,2)\oplus\u(1)\\
\su(3)\oplus\u(1)\\
   \end{aligned} $ & R \\
   \hline

$A_1+A_1+T_1$ & $\su(2)\oplus \ssl(2,\R)\oplus\u(1)$ & R\\
   \hline

\end{longtable}

\begin{longtable}{|l|c|c|}
\caption{Maximal subalgebras of $\su(2,2)$.}\label{tab:su22}
\endfirsthead
\hline
\multicolumn{3}{|l|}{\small\slshape subalgebras of $\su(2,2)$} \\
\hline
\endhead
\hline
\endfoot
\endlastfoot

\hline
type & real subalgebras & \\
\hline

$A_2$ & $\su(1,2)\oplus\u(1) (2\times)$ & R \\
\hline

$A_1+A_1+T_1$ & $
   \begin{aligned}
\ssl(2,\C)\oplus\so(1,1)\\
\ssl(2,\R)\oplus \ssl(2,\R)\oplus\u(1)\\
\su(2)\oplus \su(2)\oplus\u(1)\\
   \end{aligned} $ & R \\
   \hline

$B_2$ & $
   \begin{aligned}
\so(2,3)\\
\so(4,1)\\
   \end{aligned} $ & S \\
   \hline

$A_1+A_1$ & $
   \begin{aligned}
\su(2)\oplus \ssl(2,\R)\\
\ssl(2,\R)\oplus \ssl(2,\R)\\
   \end{aligned} $ & S \\
   \hline

\end{longtable}

\newpage

\begin{longtable}{|l|c|c|}
\caption{Maximal subalgebras of $\ssl(2,\HH)$.}\label{tab:sl2H}
\endfirsthead
\hline
\multicolumn{3}{|l|}{\small\slshape subalgebras of $\ssl(2,\HH)$} \\
\hline
\endhead
\hline
\endfoot
\endlastfoot

\hline
type & real subalgebras & \\
\hline

$B_2$ & $
   \begin{aligned}
\so(5)\\
\so(4,1)\\
   \end{aligned} $ & S \\
   \hline

$A_1+A_1+T_1$ & $
   \begin{aligned}
\ssl(2,\C)\oplus\u(1)\\
\su(2)\oplus \su(2)\oplus\so(1,1)\\
   \end{aligned} $ & R \\
   \hline

$A_1+A_1$ & $\su(2)\oplus \ssl(2,\R)$ & S\\
   \hline

\end{longtable}

\begin{longtable}{|l|c|c|}
\caption{Maximal subalgebras of $\ssl(4,\R)$.}\label{tab:sl4R}
\endfirsthead
\hline
\multicolumn{3}{|l|}{\small\slshape subalgebras of $\ssl(4,\R)$} \\
\hline
\endhead
\hline
\endfoot
\endlastfoot

\hline
type & real subalgebras & \\
\hline

$A_2$ & $\ssl(3,\R)\oplus\so(1,1)$ & R\\
\hline

$A_1+A_1+T_1$ & $
   \begin{aligned}
\ssl(2,\R)\oplus \ssl(2,\R)\oplus\so(1,1)\\
\ssl(2,\C)\oplus\u(1) (2\times)\\
   \end{aligned} $ & R \\
   \hline

$B_2$ & $\so(2,3)$ & S\\
   \hline

$A_1+A_1$ & $
   \begin{aligned}
\su(2)\oplus \su(2)\\
\ssl(2,\R)\oplus \ssl(2,\R)\\
   \end{aligned} $ & S \\
   \hline

\end{longtable}

\noindent{\bf Type $B_3$}

\begin{longtable}{|l|c|c|}
\caption{Maximal subalgebras of $\so(2,5)$.}\label{tab:so25}
\endfirsthead
\hline
\multicolumn{3}{|l|}{\small\slshape subalgebras of $\so(2,5)$} \\
\hline
\endhead
\hline
\endfoot
\endlastfoot

\hline
type & real subalgebras & \\
\hline

$A_1+A_1+A_1$ & $
   \begin{aligned}
\su(2)\oplus \ssl(2,\R)\oplus \ssl(2,\R)\\
\ssl(2,\C)\oplus \ssl(2,\R)\\
\su(2)\oplus \su(2)\oplus \ssl(2,\R)\\
   \end{aligned} $ & R \\
   \hline

$A_3$ & $
   \begin{aligned}
\su(2,2)\\
\ssl(2,\HH)\\
   \end{aligned} $ & R \\
   \hline

$B_2+T_1$ & $
   \begin{aligned}
\so(4,1)\oplus\so(1,1)\\
\so(2,3)\oplus\u(1)\\
\so(5)\oplus\u(1)\\
   \end{aligned} $ & R \\
   \hline

\end{longtable}

\begin{longtable}{|l|c|c|}
\caption{Maximal subalgebras of $\so(4,3)$.}\label{tab:so43}
\endfirsthead
\hline
\multicolumn{3}{|l|}{\small\slshape subalgebras of $\so(4,3)$} \\
\hline
\endhead
\hline
\endfoot
\endlastfoot

\hline
type & real subalgebras & \\
\hline

$A_1+A_1+A_1$ & $
   \begin{aligned}
\ssl(2,\R)\oplus \ssl(2,\R)\oplus \ssl(2,\R)\\
\ssl(2,\C)\oplus \ssl(2,\R)\\
\ssl(2,\C)\oplus \su(2)\\
\su(2)\oplus \su(2)\oplus \su(2)\\
   \end{aligned} $ & R \\
   \hline

$A_3$ & $
   \begin{aligned}
\ssl(4,\R)\\
\su(2,2)\\
   \end{aligned} $ & R \\
   \hline

$B_2+T_1$ & $
   \begin{aligned}
\so(2,3)\oplus\so(1,1)\\
\so(2,3)\oplus\u(1)\\
\so(4,1)\oplus\u(1)\\
   \end{aligned} $ & R \\
   \hline

$G_2$ & $G_{2(2)}$ & S\\
   \hline

\end{longtable}

\begin{longtable}{|l|c|c|}
\caption{Maximal subalgebras of $\so(6,1)$.}\label{tab:so61}
\endfirsthead
\hline
\multicolumn{3}{|l|}{\small\slshape subalgebras of $\so(6,1)$} \\
\hline
\endhead
\hline
\endfoot
\endlastfoot

\hline
type & real subalgebras & \\
\hline

$A_1+A_1+A_1$ & $
   \begin{aligned}
\ssl(2,\C)\oplus \su(2)\\
\su(2)\oplus \su(2)\oplus \ssl(2,\R)\\
   \end{aligned} $ & R \\
   \hline

$A_3$ & $
   \begin{aligned}
\ssl(2,\HH)\\
\su(4)\\
   \end{aligned} $ & R \\
   \hline

$B_2+T_1$ & $
   \begin{aligned}
\so(5)\oplus\so(1,1)\\
\so(4,1)\oplus\u(1)\\
   \end{aligned} $ & R \\
   \hline

\end{longtable}

\noindent{\bf Type $C_3$}

\begin{longtable}{|l|c|c|}
\caption{Maximal subalgebras of $\sp(1,2)$.}\label{tab:sp12}
\endfirsthead
\hline
\multicolumn{3}{|l|}{\small\slshape subalgebras of $\sp(1,2)$} \\
\hline
\endhead
\hline
\endfoot
\endlastfoot

\hline
type & real subalgebras & \\
\hline

$A_2$ & $\su(1,2)\oplus\u(1)$ & R\\
\hline

$A_1+B_2$ & $
   \begin{aligned}
\su(2)\oplus \so(4,1)\\
\su(2)\oplus \so(5)\\
   \end{aligned} $ & R \\
   \hline

$A_1+A_1$ & $\su(2)\oplus \ssl(2,\R)$ & S\\
   \hline

\end{longtable}

\begin{longtable}{|l|c|c|}
\caption{Maximal subalgebras of $\sp(3,\R)$.}\label{tab:sp6R}
\endfirsthead
\hline
\multicolumn{3}{|l|}{\small\slshape subalgebras of $\sp(3,\R)$} \\
\hline
\endhead
\hline
\endfoot
\endlastfoot

\hline
type & real subalgebras & \\
\hline

$A_2$ & $
   \begin{aligned}
\ssl(3,\R)\oplus\so(1,1)\\
\su(1,2)\oplus\u(1)\\
\su(3)\oplus\u(1)\\
   \end{aligned} $ & R \\
   \hline

$A_1+B_2$ & $\ssl(2,\R)\oplus \so(2,3)$ & R\\
   \hline

$A_1$ & $\ssl(2,\R)$ & S\\
   \hline

$A_1+A_1$ & $\ssl(2,\R)\oplus \ssl(2,\R)$ & S\\
   \hline

\end{longtable}

\noindent{\bf Type $A_4$}

\begin{longtable}{|l|c|c|}
\caption{Maximal subalgebras of $\su(1,4)$.}\label{tab:su14}
\endfirsthead
\hline
\multicolumn{3}{|l|}{\small\slshape subalgebras of $\su(1,4)$} \\
\hline
\endhead
\hline
\endfoot
\endlastfoot

\hline
type & real subalgebras & \\
\hline

$A_1+A_2+T_1$ & $
   \begin{aligned}
\su(2)\oplus \su(1,2)\oplus\u(1)\\
\ssl(2,\R)\oplus \su(3)\oplus\u(1)\\
   \end{aligned} $ & R \\
   \hline

$A_3+T_1$ & $
   \begin{aligned}
\su(1,3)\oplus\u(1)\\
\su(4)\oplus\u(1)\\
   \end{aligned} $ & R \\
   \hline

$B_2$ & $\so(4,1)$ & S \\
\hline
\end{longtable}

\begin{longtable}{|l|c|c|}
  \caption{Maximal subalgebras of $\su(2,3)$.}\label{tab:su23}
\endfirsthead
\hline
\multicolumn{3}{|l|}{\small\slshape subalgebras of $\su(2,3)$} \\
\hline
\endhead
\hline
\endfoot
\endlastfoot

\hline
type & real subalgebras & \\
\hline

$A_1+A_2+T_1$ & $
   \begin{aligned}
\ssl(2,\R)\oplus \su(1,2)\oplus\u(1)\\
\su(2)\oplus \su(1,2)\oplus\u(1)\\
\su(2)\oplus \su(3)\oplus\u(1)\\
   \end{aligned} $ & R \\
   \hline

$A_3+T_1$ & $
   \begin{aligned}
\su(2,2)\oplus\u(1)\\
\su(1,3)\oplus\u(1)\\
   \end{aligned} $ & R \\
   \hline

$B_2$ & $\so(2,3)$ & S \\
   \hline

\end{longtable}

\begin{longtable}{|l|c|c|}
\caption{Maximal subalgebras of $\ssl(5,\R)$.}\label{tab:sl5R}
\endfirsthead
\hline
\multicolumn{3}{|l|}{\small\slshape subalgebras of $\ssl(5,\R)$} \\
\hline
\endhead
\hline
\endfoot
\endlastfoot

\hline
type & real subalgebras & \\
\hline

$A_1+A_2+T_1$ & $\ssl(2,\R)\oplus \ssl(3,\R)\oplus\so(1,1)$ & R\\
   \hline

$A_3+T_1$ & $\ssl(4,\R)\oplus\so(1,1)$ & R \\
   \hline

$B_2$ & $
   \begin{aligned}
\so(5)\\
\so(2,3)\\
\so(4,1)\\
   \end{aligned} $ & S \\
   \hline

\end{longtable}

\noindent{\bf Type $B_4$}

\begin{longtable}{|l|c|c|}
\caption{Maximal subalgebras of $\so(2,7)$.}\label{tab:so27}
\endfirsthead
\hline
\multicolumn{3}{|l|}{\small\slshape subalgebras of $\so(2,7)$} \\
\hline
\endhead
\hline
\endfoot
\endlastfoot

\hline
type & real subalgebras & \\
\hline

$A_1+A_1+B_2$ & $
   \begin{aligned}
\ssl(2,\R)\oplus \ssl(2,\R)\oplus \so(5)\\
\su(2)\oplus \su(2)\oplus \so(2,3)\\
\ssl(2,\C)\oplus \so(4,1)\\
   \end{aligned} $ & R \\
   \hline

$A_1+A_3$ & $
   \begin{aligned}
\su(2)\oplus \su(2,2)\\
\ssl(2,\R)\oplus \ssl(2,\HH)\\
\ssl(2,\R)\oplus \su(4)\\
   \end{aligned} $ & R \\
   \hline

$B_3+T_1$ & $
   \begin{aligned}
\so(6,1)\oplus\so(1,1)\\
\so(2,5)\oplus\u(1)\\
\so(7)\oplus\u(1)\\
   \end{aligned} $ & R \\
   \hline

$D_4$ & $
   \begin{aligned}
\so^*(8)\\
\so(1,7)\\
   \end{aligned} $ & R \\
   \hline

\end{longtable}

\begin{longtable}{|l|c|c|}
\caption{Maximal subalgebras of $\so(4,5)$.}\label{tab:so45}
\endfirsthead
\hline
\multicolumn{3}{|l|}{\small\slshape subalgebras of $\so(4,5)$} \\
\hline
\endhead
\hline
\endfoot
\endlastfoot

\hline
type & real subalgebras & \\
\hline

$A_1+A_1+B_2$ & $
   \begin{aligned}
\ssl(2,\R)\oplus \ssl(2,\R)\oplus \so(2,3)\\
\ssl(2,\C)\oplus \so(2,3)\\
\ssl(2,\C)\oplus \so(4,1)\\
\su(2)\oplus \su(2)\oplus \so(4,1)\\
\su(2)\oplus \su(2)\oplus \so(5)\\
   \end{aligned} $ & R \\
   \hline

$A_1+A_3$ & $
   \begin{aligned}
\ssl(2,\R)\oplus \ssl(4,\R)\\
\ssl(2,\R)\oplus \su(2,2)\\
\su(2)\oplus \su(2,2)\\
\su(2)\oplus \ssl(2,\HH)\\
   \end{aligned} $ & R \\
   \hline

$D_4$ & $
   \begin{aligned}
\so(4,4)\\
\so(3,5)\\
   \end{aligned} $ & R \\
   \hline

$B_3+T_1$ & $
   \begin{aligned}
\so(4,3)\oplus\so(1,1)\\
\so(4,3)\oplus\u(1)\\
\so(2,5)\oplus\u(1)\\
   \end{aligned} $ & R \\
   \hline

$A_1$ & $\ssl(2,\R)$ & S\\
   \hline

$A_1+A_1$ & $\ssl(2,\R)\oplus \ssl(2,\R)$ & S\\
   \hline

\end{longtable}

\newpage

\begin{longtable}{|l|c|c|}
\caption{Maximal subalgebras of $\so(6,3)$.}\label{tab:so63}
\endfirsthead
\hline
\multicolumn{3}{|l|}{\small\slshape subalgebras of $\so(6,3)$} \\
\hline
\endhead
\hline
\endfoot
\endlastfoot

\hline
type & real subalgebras & \\
\hline

$A_1+A_1+B_2$ & $
   \begin{aligned}
\ssl(2,\R)\oplus \ssl(2,\R)\oplus \so(4,1)\\
\ssl(2,\C)\oplus \so(2,3)\\
\su(2)\oplus \su(2)\oplus \so(2,3)\\
\ssl(2,\C)\oplus \so(5)\\
   \end{aligned} $ & R \\
   \hline

$A_1+A_3$ & $
   \begin{aligned}
\su(2)\oplus \ssl(4,\R)\\
\ssl(2,\R)\oplus \su(2,2)\\
\ssl(2,\R)\oplus \ssl(2,\HH)\\
\su(2)\oplus \su(4)\\
   \end{aligned} $ & R \\
   \hline

$D_4$ & $
   \begin{aligned}
\so(3,5)\\
\so^*(8)\\
   \end{aligned} $ & R \\
   \hline

$B_3+T_1$ & $
   \begin{aligned}
\so(2,5)\oplus\so(1,1)\\
\so(4,3)\oplus\u(1)\\
\so(6,1)\oplus\u(1)\\
   \end{aligned} $ & R \\
   \hline

$A_1+A_1$ & $\su(2)\oplus\ssl(2,\R)$ & S\\
   \hline

\end{longtable}

\begin{longtable}{|l|c|c|}
\caption{Maximal subalgebras of $\so(8,1)$.}\label{tab:so81}
\endfirsthead
\hline
\multicolumn{3}{|l|}{\small\slshape subalgebras of $\so(8,1)$} \\
\hline
\endhead
\hline
\endfoot
\endlastfoot

\hline
type & real subalgebras & \\
\hline

$A_1+A_1+B_2$ & $
   \begin{aligned}
\ssl(2,\C)\oplus \so(5)\\
\su(2)\oplus \su(2)\oplus \so(4,1)\\
   \end{aligned} $ & R \\
   \hline

$A_1+A_3$ & $
   \begin{aligned}
\su(2)\oplus \ssl(2,\HH)\\
\ssl(2,\R)\oplus \su(4)\\
   \end{aligned} $ & R \\
   \hline

$D_4$ & $
   \begin{aligned}
\so(1,7)\\
\so(8)\\
   \end{aligned} $ & R \\
   \hline

$B_3+T_1$ & $
   \begin{aligned}
\so(7)\oplus\so(1,1)\\
\so(6,1)\oplus\u(1)\\
   \end{aligned} $ & R \\
   \hline

\end{longtable}

\noindent{\bf Type $C_4$}

\begin{longtable}{|l|c|c|}
\caption{Maximal subalgebras of $\sp(1,3)$.}\label{tab:sp13}
\endfirsthead
\hline
\multicolumn{3}{|l|}{\small\slshape subalgebras of $\sp(1,3)$} \\
\hline
\endhead
\hline
\endfoot
\endlastfoot

\hline
type & real subalgebras & \\
\hline

$A_3+T_1$ & $\su(1,3)\oplus\u(1)$ & R \\
   \hline

$B_2+B_2$ & $\so(5)\oplus \so(4,1)$ & R\\
   \hline

$A_1+C_3$ & $
   \begin{aligned}
\su(2)\oplus \sp(1,2)\\
\su(2)\oplus \sp(3)\\
   \end{aligned} $ & R \\
   \hline

\end{longtable}

\begin{longtable}{|l|c|c|}
\caption{Maximal subalgebras of $\sp(2,2)$.}\label{tab:sp22}
\endfirsthead
\hline
\multicolumn{3}{|l|}{\small\slshape subalgebras of $\sp(2,2)$} \\
\hline
\endhead
\hline
\endfoot
\endlastfoot

\hline
type & real subalgebras & \\
\hline

$A_3+T_1$ & $
   \begin{aligned}
\ssl(2,\HH)\oplus\so(1,1)\\
\su(2,2)\oplus\u(1)\\
   \end{aligned} $ & R \\
   \hline

$B_2+B_2$ & $
   \begin{aligned}
\so(4,1)\oplus \so(4,1)\\
\so(5,\C)\\
\so(5)\oplus \so(5)\\
   \end{aligned} $ & R \\
   \hline

$A_1+C_3$ & $\su(2)\oplus \sp(1,2) (2\times)$ & R \\
   \hline

$A_1+A_1+A_1$ & $\su(2)\oplus \ssl(2,\R)\oplus \ssl(2,\R)$ & S \\
   \hline

\end{longtable}

\begin{longtable}{|l|c|c|}
\caption{Maximal subalgebras of $\sp(4,\R)$.}\label{tab:sp8R}
\endfirsthead
\hline
\multicolumn{3}{|l|}{\small\slshape subalgebras of $\sp(4,\R)$} \\
\hline
\endhead
\hline
\endfoot
\endlastfoot

\hline
type & real subalgebras & \\
\hline

$A_3+T_1$ & $
   \begin{aligned}
\ssl(4,\R)\oplus\so(1,1)\\
\su(2,2)\oplus\u(1)\\
\su(1,3)\oplus\u(1)\\
\su(4)\oplus\u(1)\\
   \end{aligned} $ & R \\
   \hline

$B_2+B_2$ & $
   \begin{aligned}
\so(2,3)\oplus \so(2,3)\\
\so(5,\C)\\
   \end{aligned} $ & R \\
   \hline

$A_1+C_3$ & $\ssl(2,\R)\oplus \sp(3,\R)$ & R \\
   \hline

$A_1$ & $\ssl(2,\R)$ & S\\
   \hline

$A_1+A_1+A_1$ & $
   \begin{aligned}
\su(2)\oplus \su(2)\oplus \ssl(2,\R)\\
\ssl(2,\R)\oplus \ssl(2,\R)\oplus \ssl(2,\R)\\
   \end{aligned} $ & S \\
   \hline

\end{longtable}

\noindent{\bf Type $D_4$}

\begin{longtable}{|l|c|c|}
\caption{Maximal subalgebras of $\so^*(8)$.}\label{tab:so*8}
\endfirsthead
\hline
\multicolumn{3}{|l|}{\small\slshape subalgebras of $\so^*(8)$} \\
\hline
\endhead
\hline
\endfoot
\endlastfoot

\hline
type & real subalgebras & \\
\hline

$A_1+A_1+A_1+A_1$ & $
   \begin{aligned}
\su(2)\oplus \su(2)\oplus \ssl(2,\R)\oplus \ssl(2,\R)\\
\ssl(2,\C)\oplus \ssl(2,\C)\\
   \end{aligned} $ & R \\
   \hline

$A_3+T_1$ & $
   \begin{aligned}
\su(2,2)\oplus\u(1)\\
\ssl(2,\HH)\oplus\so(1,1)\\
\su(4)\oplus\u(1)\\
   \end{aligned} $ & R \\
   \hline

$A_2$ & $
   \begin{aligned}
\su(2)\oplus \so(2,3)\\
\ssl(2,\R)\oplus \so(5)\\
\ssl(2,\R)\oplus \so(4,1)\\
   \end{aligned} $ & S \\
   \hline

$B_3$ & $
   \begin{aligned}
\so(2,5)\\
\so(6,1)\\
   \end{aligned} $ & S \\
   \hline

\end{longtable}

\begin{longtable}{|l|c|c|}
\caption{Maximal subalgebras of $\so(4,4)$.}\label{tab:so44}
\endfirsthead
\hline
\multicolumn{3}{|l|}{\small\slshape subalgebras of $\so(4,4)$} \\
\hline
\endhead
\hline
\endfoot
\endlastfoot

\hline
type & real subalgebras & \\
\hline

$A_1+A_1+A_1+A_1$ & $
   \begin{aligned}
\ssl(2,\R)\oplus \ssl(2,\R)\oplus \ssl(2,\R)\oplus \ssl(2,\R)\\
\ssl(2,\C)\oplus \ssl(2,\C) (3\times)\\
\su(2)\oplus \su(2)\oplus \su(2)\oplus \su(2)\\
   \end{aligned} $ & R \\
   \hline

$A_3+T_1$ & $
   \begin{aligned}
\ssl(4,\R)\oplus\so(1,1)\\
\su(2,2)\oplus\u(1) (6\times)\\
   \end{aligned} $ & R \\
   \hline

$A_2$ & $\su(1,2)$ & S\\
   \hline

$A_1+B_2$ & $
   \begin{aligned}
\su(2)\oplus \so(4,1)\\
\ssl(2,\R)\oplus \so(2,3)\\
   \end{aligned} $ & S \\
   \hline

$B_3$ & $\so(4,3)$ & S\\
   \hline

\end{longtable}

\begin{longtable}{|l|c|c|}
\caption{Maximal subalgebras of $\so(3,5)$.}\label{tab:so35}
\endfirsthead
\hline
\multicolumn{3}{|l|}{\small\slshape subalgebras of $\so(3,5)$} \\
\hline
\endhead
\hline
\endfoot
\endlastfoot

\hline
type & real subalgebras & \\
\hline

$A_1+A_1+A_1+A_1$ & $
   \begin{aligned}
\ssl(2,\C)\oplus \ssl(2,\R)\oplus \ssl(2,\R)\\
\ssl(2,\C)\oplus \su(2)\oplus \su(2)\\
   \end{aligned} $ & R \\
   \hline

$A_3+T_1$ & $
   \begin{aligned}
\ssl(4,\R)\oplus\u(1)\\
\su(2,2)\oplus\so(1,1)\\
\ssl(2,\HH)\oplus\u(1)\\
   \end{aligned} $ & R \\
   \hline

$A_2$ & $\ssl(3,\R)$ & S\\
   \hline

$A_1+B_2$ & $
   \begin{aligned}
\su(2)\oplus \so(5)\\
\su(2)\oplus \so(2,3)\\
\ssl(2,\R)\oplus \so(2,3)\\
\ssl(2,\R)\oplus \so(4,1)\\
   \end{aligned} $ & S \\
   \hline

$B_3$ & $
   \begin{aligned}
\so(2,5)\\
\so(4,3)\\
   \end{aligned} $ & S \\
   \hline

\end{longtable}

\begin{longtable}{|l|c|c|}
\caption{Maximal subalgebras of $\so(1,7)$.}\label{tab:so17}
\endfirsthead
\hline
\multicolumn{3}{|l|}{\small\slshape subalgebras of $\so(1,7)$} \\
\hline
\endhead
\hline
\endfoot
\endlastfoot

\hline
type & real subalgebras & \\
\hline

$A_1+A_1+A_1+A_1$ & $\ssl(2,\C)\oplus \su(2)\oplus \su(2)$ & R\\
\hline

$A_3+T_1$ & $
   \begin{aligned}
\ssl(2,\HH)\oplus\u(1)\\
\su(4)\oplus\so(1,1)\\
   \end{aligned} $ & R \\
   \hline

$A_1+B_2$ & $
   \begin{aligned}
\su(2)\oplus \so(4,1)\\
\ssl(2,\R)\oplus \so(5)\\
   \end{aligned} $ & S \\
   \hline

$B_3$ & $
   \begin{aligned}
\so(7)\\
\so(6,1)\\
   \end{aligned} $ & S \\
   \hline

\end{longtable}

\noindent{\bf Type $F_4$}

\begin{longtable}{|l|c|c|}
\caption{Maximal subalgebras of $F_{4(4)}$.}\label{tab:F44}
\endfirsthead
\hline
\multicolumn{3}{|l|}{\small\slshape subalgebras of $F_{4(4)}$} \\
\hline
\endhead
\hline
\endfoot
\endlastfoot

\hline
type & real subalgebras & \\
\hline

$A_2+A_2$ & $
   \begin{aligned}
     \ssl(3,\R)\oplus \ssl(3,\R)\\
     \su(1,2)\oplus\su(1,2)\\
     \su(3)\oplus \su(1,2)\\
   \end{aligned} $ & R \\
   \hline

$A_1+C_3$ & $
   \begin{aligned}
\ssl(2,\R)\oplus \sp(3,\R)\\
\su(2)\oplus \sp(1,2)\\
\su(2)\oplus \sp(3)\\
   \end{aligned} $ & R \\
   \hline

$B_4$ & $\so(4,5)$ & R\\
   \hline

$A_1$ & $\ssl(2,\R)$ & S\\
   \hline

$A_1+G_2$ & $
   \begin{aligned}
     \su(2)\oplus G_{2(2)}\\
     \ssl(2,\R)\oplus G_{2(2)}\\
   \end{aligned} $ & S \\
   \hline

\end{longtable}

\begin{longtable}{|l|c|c|}
\caption{Maximal subalgebras of $F_{4(-20)}$.}\label{tab:F420}
\endfirsthead
\hline
\multicolumn{3}{|l|}{\small\slshape subalgebras of $F_{4(-20)}$} \\
\hline
\endhead
\hline
\endfoot
\endlastfoot

\hline
type & real subalgebras & \\
\hline

$A_2+A_2$ & $\su(3)\oplus \su(1,2)$ & R\\
   \hline

$A_1+C_3$ & $\su(2)\oplus \sp(1,2)$ & R\\
   \hline

$B_4$ & $
   \begin{aligned}
\so(8,1)\\
\so(9)\\
   \end{aligned} $ & R \\
   \hline

$A_1+G_2$ & $\ssl(2,\R)\oplus G_2^c$ & S\\
   \hline

\end{longtable}

\noindent{\bf Type $A_5$}

\begin{longtable}{|l|c|c|}
\caption{Maximal subalgebras of $\su(1,5)$.}\label{tab:su15}
\endfirsthead
\hline
\multicolumn{3}{|l|}{\small\slshape subalgebras of $\su(1,5)$} \\
\hline
\endhead
\hline
\endfoot
\endlastfoot

\hline
type & real subalgebras & \\
\hline

$A_2+A_2+T_1$ & $\su(3)\oplus \su(1,2)\oplus\u(1)$ & R \\
   \hline

$A_1+A_3+T_1$ & $
   \begin{aligned}
\su(2)\oplus \su(1,3)\oplus\u(1)\\
\ssl(2,\R)\oplus \su(4)\oplus\u(1)\\
   \end{aligned} $ & R \\
   \hline

$A_4+T_1$ & $
   \begin{aligned}
\su(1,4)\oplus\u(1)\\
\su(5)\oplus\u(1)\\
   \end{aligned} $ & R \\
   \hline

$A_3$ & $\ssl(2,\HH)$ & S \\
   \hline

\end{longtable}

\begin{longtable}{|l|c|c|}
\caption{Maximal subalgebras of $\su(2,4)$.}\label{tab:su24}
\endfirsthead
\hline
\multicolumn{3}{|l|}{\small\slshape subalgebras of $\su(2,4)$} \\
\hline
\endhead
\hline
\endfoot
\endlastfoot

\hline
type & real subalgebras & \\
\hline

$A_2+A_2+T_1$ & $
   \begin{aligned}
\su(1,2)\oplus \su(1,2)\oplus\u(1)\\
\su(3)\oplus \su(1,2)\oplus\u(1)\\
   \end{aligned} $ & R \\
   \hline

$A_1+A_3+T_1$ & $
   \begin{aligned}
\su(2)\oplus \su(2,2)\oplus\u(1)\\
\ssl(2,\R)\oplus \su(1,3)\oplus\u(1)\\
\su(2)\oplus \su(4)\oplus\u(1)\\
   \end{aligned} $ & R \\
   \hline

$A_4+T_1$ & $
   \begin{aligned}
\su(2,3)\oplus\u(1)\\
\su(1,4)\oplus\u(1)\\
   \end{aligned} $ & R \\
   \hline

$A_2$ & $\su(1,2)$ & S \\
   \hline

$A_1+A_2$ & $\su(2)\oplus \su(1,2)$ & S \\
   \hline

$A_3$ & $\su(2,2)$ & S \\
   \hline

$C_3$ & $\sp(1,2)$ & S \\
   \hline

\end{longtable}

\begin{longtable}{|l|c|c|}
\caption{Maximal subalgebras of $\su(3,3)$.}\label{tab:su33}
\endfirsthead
\hline
\multicolumn{3}{|l|}{\small\slshape subalgebras of $\su(3,3)$} \\
\hline
\endhead
\hline
\endfoot
\endlastfoot

\hline
type & real subalgebras & \\
\hline

$A_2+A_2+T_1$ & $
   \begin{aligned}
\ssl(3,\C)\oplus\so(1,1)\\
\su(1,2)\oplus \su(1,2)\oplus\u(1)\\
\su(3)\oplus \su(3)\oplus\u(1)\\
   \end{aligned} $ & R \\
   \hline

$A_1+A_3+T_1$ & $
   \begin{aligned}
\ssl(2,\R)\oplus \su(2,2)\oplus\u(1)\\
\su(2)\oplus \su(1,3)\oplus\u(1) (2\times)\\
   \end{aligned} $ & R \\
   \hline

$A_4+T_1$ & $\su(2,3)\oplus\u(1) (2\times)$ & R \\
   \hline

$A_1+A_2$ & $
   \begin{aligned}
\ssl(2,\R)\oplus \su(3)\\
\ssl(2,\R)\oplus \su(1,2)\\
   \end{aligned} $ & S \\
   \hline

$A_3$ & $
   \begin{aligned}
\su(1,3)\\
\ssl(4,\R)\\
   \end{aligned} $ & S \\
   \hline

$C_3$ & $\sp(3,\R)$ & S \\
   \hline

\end{longtable}

\newpage

\begin{longtable}{|l|c|c|}
\caption{Maximal subalgebras of $\ssl(3,\HH)$.}\label{tab:sl3H}
\endfirsthead
\hline
\multicolumn{3}{|l|}{\small\slshape subalgebras of $\ssl(3,\HH)$} \\
\hline
\endhead
\hline
\endfoot
\endlastfoot

\hline
type & real subalgebras & \\
\hline

$A_2+A_2+T_1$ & $\ssl(3,\C)\oplus\u(1)$ & S \\
   \hline

$A_1+A_3+T_1$ & $\su(2)\oplus \ssl(2,\HH)\oplus\so(1,1)$ & R\\
   \hline

$A_1+A_2$ & $\su(2)\oplus \ssl(3,\R)$ & S\\
   \hline

$A_3$ & $\su(1,3)$ & S\\
   \hline

$C_3$ & $
   \begin{aligned}
\sp(3)\\
\sp(1,2)\\
   \end{aligned} $ & S \\
   \hline

\end{longtable}

\begin{longtable}{|l|c|c|}
\caption{Maximal subalgebras of $\ssl(6,\R)$.}\label{tab:sl6R}
\endfirsthead
\hline
\multicolumn{3}{|l|}{\small\slshape subalgebras of $\ssl(6,\R)$} \\
\hline
\endhead
\hline
\endfoot
\endlastfoot

\hline
type & real subalgebras & \\
\hline

$A_2+A_2+T_1$ & $
   \begin{aligned}
\ssl(3,\R)\oplus \ssl(3,\R)\oplus\so(1,1)\\
\ssl(3,\C)\oplus\u(1)\\
   \end{aligned} $ & R \\
   \hline

$A_1+A_3+T_1$ & $\ssl(2,\R)\oplus \ssl(4,\R)\oplus\so(1,1)$ & R\\
   \hline

$A_4+T_1$ & $\ssl(5,\R)\oplus\so(1,1)$ & R\\
   \hline

$A_2$ & $\ssl(3,\R)$ & S \\
   \hline

$A_1+A_2$ & $\ssl(2,\R)\oplus \ssl(3,\R)$ & S \\
   \hline

$A_3$ & $
   \begin{aligned}
\su(4)\\
\su(2,2)\\
\ssl(2,\HH)\\
\ssl(4,\R)\\
   \end{aligned} $ & S \\
   \hline

$C_3$ & $\sp(3,\R)$ & S\\
   \hline

\end{longtable}

\noindent{\bf Type $B_5$}

\begin{longtable}{|l|c|c|}
\caption{Maximal subalgebras of $\so(2,9)$.}\label{tab:so29}
\endfirsthead
\hline
\multicolumn{3}{|l|}{\small\slshape subalgebras of $\so(2,9)$} \\
\hline
\endhead
\hline
\endfoot
\endlastfoot

\hline
type & real subalgebras & \\
\hline

$A_3+B_2$ & $
   \begin{aligned}
\su(2,2)\oplus \so(5)\\
\ssl(2,\HH)\oplus \so(4,1)\\
\su(4)\oplus \so(2,3)\\
   \end{aligned} $ & R \\
   \hline

$A_1+A_1+B_3$ & $
   \begin{aligned}
\ssl(2,\R)\oplus \ssl(2,\R)\oplus \so(7)\\
\ssl(2,\C)\oplus \so(6,1)\\
\su(2)\oplus \su(2)\oplus \so(2,5)\\
   \end{aligned} $ & R \\
   \hline

$A_1+D_4$ & $
   \begin{aligned}
\su(2)\oplus \so^*(8)\\
\ssl(2,\R)\oplus \so(1,7)\\
\ssl(2,\R)\oplus \so(8)\\
   \end{aligned} $ & R \\
   \hline

$D_5$ & $
   \begin{aligned}
\so(2,8)\\
\so(9,1)\\
   \end{aligned} $ & R \\
   \hline

$B_4+T_1$ & $
   \begin{aligned}
\so(8,1)\oplus\so(1,1)\\
\so(2,7)\oplus\u(1)\\
\so(9)\oplus\u(1)\\
   \end{aligned} $ & R \\
   \hline

\end{longtable}

\newpage

\begin{longtable}{|l|c|c|}
\caption{Maximal subalgebras of $\so(4,7)$.}\label{tab:so47}
\endfirsthead
\hline
\multicolumn{3}{|l|}{\small\slshape subalgebras of $\so(4,7)$} \\
\hline
\endhead
\hline
\endfoot
\endlastfoot

\hline
type & real subalgebras & \\
\hline

$A_3+B_2$ & $
   \begin{aligned}
\ssl(4,\R)\oplus \so(4,1)\\
\su(2,2)\oplus \so(2,3)\\
\ssl(2,\HH)\oplus \so(2,3)\\
\su(2,2)\oplus \so(5)\\
\su(4)\oplus \so(4,1)\\
   \end{aligned} $ & R \\
   \hline

$A_1+A_1+B_3$ & $
   \begin{aligned}
\ssl(2,\R)\oplus \ssl(2,\R)\oplus \so(2,5)\\
\ssl(2,\C)\oplus \so(4,3)\\
\su(2)\oplus \su(2)\oplus \so(4,3)\\
\ssl(2,\C)\oplus \so(6,1)\\
\su(2)\oplus \su(2)\oplus \so(7)\\
   \end{aligned} $ & R \\
   \hline

$A_1+D_4$ & $
   \begin{aligned}
\su(2)\oplus \so(4,4)\\
\ssl(2,\R)\oplus \so(3,5)\\
\ssl(2,\R)\oplus \so^*(8)\\
\su(2)\oplus \so(1,7)\\
   \end{aligned} $ & R \\
   \hline

$D_5$ & $
   \begin{aligned}
\so(4,6)\\
\so(3,7)\\
   \end{aligned} $ & R \\
   \hline

$B_4+T_1$ & $
   \begin{aligned}
\so(6,3)\oplus\so(1,1)\\
\so(4,5)\oplus\u(1)\\
\so(2,7)\oplus\u(1)\\
   \end{aligned} $ & R \\
   \hline

\end{longtable}

\begin{longtable}{|l|c|c|}
\caption{Maximal subalgebras of $\so(6,5)$.}\label{tab:so65}
\endfirsthead
\hline
\multicolumn{3}{|l|}{\small\slshape subalgebras of $\so(6,5)$} \\
\hline
\endhead
\hline
\endfoot
\endlastfoot

\hline
type & real subalgebras & \\
\hline

$A_3+B_2$ & $
   \begin{aligned}
\ssl(4,\R)\oplus \so(2,3)\\
\su(2,2)\oplus \so(2,3)\\
\su(2,2)\oplus \so(4,1)\\
\ssl(2,\HH)\oplus \so(4,1)\\
\ssl(2,\HH)\oplus \so(5)\\
\su(4)\oplus \so(5)\\
   \end{aligned} $ & R \\
   \hline

$A_1+A_1+B_3$ & $
   \begin{aligned}
\ssl(2,\R)\oplus \ssl(2,\R)\oplus \so(4,3)\\
\ssl(2,\C)\oplus \so(4,3)\\
\ssl(2,\C)\oplus \so(2,5)\\
\su(2)\oplus \su(2)\oplus \so(2,5)\\
\su(2)\oplus \su(2)\oplus \so(6,1)\\
   \end{aligned} $ & R \\
   \hline

$A_1+D_4$ & $
   \begin{aligned}
\ssl(2,\R)\oplus \so(4,4)\\
\ssl(2,\R)\oplus \so(3,5)\\
\su(2)\oplus \so(3,5)\\
\su(2)\oplus \so^*(8)\\
   \end{aligned} $ & R \\
   \hline

$D_5$ & $
   \begin{aligned}
\so(5,5)\\
\so(4,6)\\
   \end{aligned} $ & R \\
   \hline

$B_4+T_1$ & $
   \begin{aligned}
\so(4,5)\oplus\so(1,1)\\
\so(4,5)\oplus\u(1)\\
\so(6,3)\oplus\u(1)\\
   \end{aligned} $ & R \\
   \hline

$A_1$ & $\ssl(2,\R)$ & S\\
   \hline

\end{longtable}

\begin{longtable}{|l|c|c|}
\caption{Maximal subalgebras of $\so(8,3)$.}\label{tab:so83}
\endfirsthead
\hline
\multicolumn{3}{|l|}{\small\slshape subalgebras of $\so(8,3)$} \\
\hline
\endhead
\hline
\endfoot
\endlastfoot

\hline
type & real subalgebras & \\
\hline

$A_3+B_2$ & $
   \begin{aligned}
\ssl(4,\R)\oplus \so(5)\\
\su(2,2)\oplus \so(4,1)\\
\ssl(2,\HH)\oplus \so(2,3)\\
\su(4)\oplus \so(2,3)\\
   \end{aligned} $ & R \\
   \hline

$A_1+A_1+B_3$ & $
   \begin{aligned}
\ssl(2,\R)\oplus \ssl(2,\R)\oplus \so(6,1)\\
\ssl(2,\C)\oplus \so(2,5)\\
\su(2)\oplus \su(2)\oplus \so(4,3)\\
\ssl(2,\C)\oplus \so(7)\\
   \end{aligned} $ & R \\
   \hline

$A_1+D_4$ & $
   \begin{aligned}
\su(2)\oplus \so(3,5)\\
\ssl(2,\R)\oplus \so^*(8)\\
\ssl(2,\R)\oplus \so(1,7)\\
\su(2)\oplus \so(8)\\
   \end{aligned} $ & R \\
   \hline

$D_5$ & $
   \begin{aligned}
\so(3,7)\\
\so(2,8)\\
   \end{aligned} $ & R \\
   \hline

$B_4+T_1$ & $
   \begin{aligned}
\so(2,7)\oplus\so(1,1)\\
\so(6,3)\oplus\u(1)\\
\so(8,1)\oplus\u(1)\\
   \end{aligned} $ & R \\
   \hline

\end{longtable}

\begin{longtable}{|l|c|c|}
\caption{Maximal subalgebras of $\so(10,1)$.}\label{tab:so101}
\endfirsthead
\hline
\multicolumn{3}{|l|}{\small\slshape subalgebras of $\so(10,1)$} \\
\hline
\endhead
\hline
\endfoot
\endlastfoot

\hline
type & real subalgebras & \\
\hline

$A_3+B_2$ & $
   \begin{aligned}
\ssl(2,\HH)\oplus \so(5)\\
\su(4)\oplus \so(4,1)\\
   \end{aligned} $ & R \\
   \hline

$A_1+A_1+B_3$ & $
   \begin{aligned}
\ssl(2,\C)\oplus \so(7)\\
\su(2)\oplus \su(2)\oplus \so(6,1)\\
   \end{aligned} $ & R \\
   \hline

$A_1+D_4$ & $
   \begin{aligned}
\su(2)\oplus \so(1,7)\\
\ssl(2,\R)\oplus \so(8)\\
   \end{aligned} $ & R \\
   \hline

$D_5$ & $
   \begin{aligned}
\so(9,1)\\
\so(10)\\
   \end{aligned} $ & R \\
   \hline

$B_4+T_1$ & $
   \begin{aligned}
\so(9)\oplus\so(1,1)\\
\so(8,1)\oplus\u(1)\\
   \end{aligned} $ & R \\
   \hline

\end{longtable}

\noindent{\bf Type $C_5$}

\begin{longtable}{|l|c|c|}
\caption{Maximal subalgebras of $\sp(1,4)$.}\label{tab:sp14}
\endfirsthead
\hline
\multicolumn{3}{|l|}{\small\slshape subalgebras of $\sp(1,4)$} \\
\hline
\endhead
\hline
\endfoot
\endlastfoot

\hline
type & real subalgebras & \\
\hline

$A_4+T_1$ & $\su(1,4)\oplus\u(1)$ & R\\
   \hline

$B_2+C_3$ & $
   \begin{aligned}
\so(4,1)\oplus \sp(3)\\
\so(5)\oplus \sp(1,2)\\
   \end{aligned} $ & R \\
   \hline

$A_1+C_4$ & $
   \begin{aligned}
\su(2)\oplus \sp(1,3)\\
\su(2)\oplus \sp(4)\\
   \end{aligned} $ & R \\
   \hline

$A_1+B_2$ & $\su(2)\oplus \so(4,1)$ & S \\
   \hline

\end{longtable}

\begin{longtable}{|l|c|c|}
\caption{Maximal subalgebras of $\sp(2,3)$.}\label{tab:sp23}
\endfirsthead
\hline
\multicolumn{3}{|l|}{\small\slshape subalgebras of $\sp(2,3)$} \\
\hline
\endhead
\hline
\endfoot
\endlastfoot

\hline
type & real subalgebras & \\
\hline

$A_4+T_1$ & $\su(2,3)$ & R \\
   \hline

$B_2+C_3$ & $
   \begin{aligned}
\so(4,1)\oplus \sp(1,2)\\
\so(5)\oplus \sp(1,2)\\
\so(5)\oplus \sp(3)\\
   \end{aligned} $ & R \\
   \hline

$A_1+C_4$ & $
   \begin{aligned}
\su(2)\oplus \sp(2,2)\\
\su(2)\oplus \sp(1,3)\\
   \end{aligned} $ & R \\
   \hline

$A_1+B_2$ & $\su(2)\oplus \so(2,3)$ & S\\
   \hline

\end{longtable}

\begin{longtable}{|l|c|c|}
\caption{Maximal subalgebras of $\sp(5,\R)$.}\label{tab:sp10R}
\endfirsthead
\hline
\multicolumn{3}{|l|}{\small\slshape subalgebras of $\sp(5,\R)$} \\
\hline
\endhead
\hline
\endfoot
\endlastfoot

\hline
type & real subalgebras & \\
\hline

$A_4+T_1$ & $
   \begin{aligned}
\ssl(5,\R)\oplus\so(1,1)\\
\su(2,3)\oplus\u(1)\\
\su(1,4)\oplus\u(1)\\
\su(5)\oplus\u(1)\\
   \end{aligned} $ & R \\
   \hline

$B_2+C_3$ & $\so(2,3)\oplus \sp(3,\R)$ & R \\
   \hline

$A_1+C_4$ & $\ssl(2,\R)\oplus \sp(4,\R)$ & R\\
   \hline

$A_1+B_2$ & $
   \begin{aligned}
     \ssl(2,\R)\oplus\so(5)\\
     \ssl(2,\R)\oplus \so(2,3)\\
\ssl(2,\R)\oplus \so(4,1)\\
   \end{aligned} $ & S \\
   \hline

$A_1$ & $\ssl(2,\R)$ & S\\
\hline
\end{longtable}

\noindent{\bf Type $D_5$}

\begin{longtable}{|l|c|c|}
\caption{Maximal subalgebras of $\so(2,8)$.}\label{tab:so28}
\endfirsthead
\hline
\multicolumn{3}{|l|}{\small\slshape subalgebras of $\so(2,8)$} \\
\hline
\endhead
\hline
\endfoot
\endlastfoot

\hline
type & real subalgebras & \\
\hline

$A_4+T_1$ & $\su(1,4)\oplus\u(1) (2\times)$ & R \\
   \hline

$A_1+A_1+A_3$ & $
   \begin{aligned}
\ssl(2,\R)\oplus \ssl(2,\R)\oplus \su(4)\\
\ssl(2,\C)\oplus \ssl(2,\HH)\\
\su(2)\oplus \su(2)\oplus \su(2,2)\\
   \end{aligned} $ & R \\
   \hline

$D_4+T_1$ & $
   \begin{aligned}
\so(1,7)\oplus\so(1,1)\\
\so^*(8)\oplus\u(1)\\
\so(8)\oplus\u(1)\\
   \end{aligned} $ & R \\
   \hline

$B_2+B_2$ & $
   \begin{aligned}
\so(5)\oplus \so(2,3)\\
\so(4,1)\oplus \so(4,1)\\
   \end{aligned} $ & S \\
   \hline

$A_1+B_3$ & $
   \begin{aligned}
\su(2)\oplus \so(2,5)\\
\ssl(2,\R)\oplus \so(7)\\
\ssl(2,\R)\oplus \so(6,1)\\
   \end{aligned} $ & S \\
   \hline

$B_4$ & $
   \begin{aligned}
\so(2,7)\\
\so(8,1)\\
   \end{aligned} $ & S \\
   \hline

\end{longtable}

\begin{longtable}{|l|c|c|}
\caption{Maximal subalgebras of $\so(4,6)$.}\label{tab:so46}
\endfirsthead
\hline
\multicolumn{3}{|l|}{\small\slshape subalgebras of $\so(4,6)$} \\
\hline
\endhead
\hline
\endfoot
\endlastfoot

\hline
type & real subalgebras & \\
\hline

$A_4+T_1$ & $\su(2,3)\oplus\u(1) (2\times)$ & R \\
   \hline

$A_1+A_1+A_3$ & $
   \begin{aligned}
\ssl(2,\R)\oplus \ssl(2,\R)\oplus \su(2,2)\\
\ssl(2,\C)\oplus \ssl(4,\R)\\
\ssl(2,\C)\oplus \ssl(2,\HH)\\
\su(2)\oplus \su(2)\oplus \su(2,2)\\
\su(2)\oplus \su(2)\oplus \su(4)\\
   \end{aligned} $ & R \\
   \hline

$D_4+T_1$ & $
   \begin{aligned}
\so(3,5)\oplus\so(1,1)\\
\so(4,4)\oplus\u(1)\\
\so^*(8)\oplus\u(1)\\
   \end{aligned} $ & R \\
   \hline

$B_2$ & $
   \begin{aligned}
\so(2,3)\\
\so(4,1)\\
   \end{aligned} $ & S \\
   \hline

$B_2+B_2$ & $
   \begin{aligned}
\so(5)\oplus \so(4,1)\\
\so(2,3)\oplus \so(2,3)\\
\so(2,3)\oplus \so(4,1)\\
   \end{aligned} $ & S \\
   \hline

$A_1+B_3$ & $
   \begin{aligned}
\su(2)\oplus \so(4,3)\\
\su(2)\oplus \so(6,1)\\
\ssl(2,\R)\oplus \so(2,5)\\
\ssl(2,\R)\oplus \so(4,3)\\
   \end{aligned} $ & S \\
   \hline

$B_4$ & $
   \begin{aligned}
\so(4,5)\\
\so(6,3)\\
   \end{aligned} $ & S \\
   \hline

\end{longtable}

\begin{longtable}{|l|c|c|}
\caption{Maximal subalgebras of $\so^*(10)$.}\label{tab:so*10}
\endfirsthead
\hline
\multicolumn{3}{|l|}{\small\slshape subalgebras of $\so^*(10)$} \\
\hline
\endhead
\hline
\endfoot
\endlastfoot

\hline
type & real subalgebras & \\
\hline

$A_4+T_1$ & $
   \begin{aligned}
\su(2,3)\oplus\u(1)\\
\su(1,4)\oplus\u(1)\\
\su(5)\oplus\u(1)\\
   \end{aligned} $ & R \\
   \hline

$A_1+A_1+A_3$ & $\su(2)\oplus \ssl(2,\R)\oplus \su(1,3)$ & R \\
   \hline

$D_4+T_1$ & $\so^*(8)\oplus\u(1)$ & R \\
   \hline

\end{longtable}

\begin{longtable}{|l|c|c|}
\caption{Maximal subalgebras of $\so(9,1)$.}\label{tab:so91}
\endfirsthead
\hline
\multicolumn{3}{|l|}{\small\slshape subalgebras of $\so(9,1)$} \\
\hline
\endhead
\hline
\endfoot
\endlastfoot

\hline
type & real subalgebras & \\
\hline

$A_1+A_1+A_3$ & $
   \begin{aligned}
\su(2)\oplus \su(2)\oplus \ssl(2,\HH)\\
\ssl(2,\C)\oplus \su(4)\\
   \end{aligned} $ & R \\
   \hline

$D_4+T_1$ & $
   \begin{aligned}
\so(1,7)\oplus\u(1)\\
\so(8)\oplus\so(1,1)\\
   \end{aligned} $ & R \\
   \hline

$B_2+B_2$ & $\so(5)\oplus \so(4,1)$ & S\\
   \hline

$A_1+B_3$ & $
   \begin{aligned}
\su(2)\oplus \so(6,1)\\
\ssl(2,\R)\oplus \so(7)\\
   \end{aligned} $ & S \\
   \hline

$B_4$ & $
   \begin{aligned}
\so(9)\\
\so(8,1)\\
   \end{aligned} $ & S \\
   \hline

\end{longtable}

\begin{longtable}{|l|c|c|}
\caption{Maximal subalgebras of $\so(3,7)$.}\label{tab:so37}
\endfirsthead
\hline
\multicolumn{3}{|l|}{\small\slshape subalgebras of $\so(3,7)$} \\
\hline
\endhead
\hline
\endfoot
\endlastfoot

\hline
type & real subalgebras & \\
\hline

$A_1+A_1+A_3$ & $
   \begin{aligned}
\ssl(2,\R)\oplus \ssl(2,\R)\oplus \ssl(2,\HH)\\
\ssl(2,\C)\oplus \su(2,2)\\
\su(2)\oplus \su(2)\oplus \ssl(4,\R)\\
\ssl(2,\C)\oplus \su(4)\\
   \end{aligned} $ & R \\
   \hline

$D_4+T_1$ & $
   \begin{aligned}
\so^*(8)\oplus\so(1,1)\\
\so(3,5)\oplus\u(1)\\
\so(1,7)\oplus\u(1)\\
   \end{aligned} $ & R \\
   \hline

$B_2+B_2$ & $
   \begin{aligned}
\so(5)\oplus \so(2,3)\\
\so(2,3)\oplus \so(4,1)\\
   \end{aligned} $ & S \\
   \hline

$A_1+B_3$ & $
   \begin{aligned}
\su(2)\oplus \so(7)\\
\su(2)\oplus \so(4,3)\\
\ssl(2,\R)\oplus \so(2,5)\\
\ssl(2,\R)\oplus \so(6,1)\\
   \end{aligned} $ & S \\
   \hline

$B_4$ & $
   \begin{aligned}
\so(2,7)\\
\so(6,3)\\
   \end{aligned} $ & S \\
   \hline

\end{longtable}

\begin{longtable}{|l|c|c|}
\caption{Maximal subalgebras of $\so(5,5)$.}\label{tab:so55}
\endfirsthead
\hline
\multicolumn{3}{|l|}{\small\slshape subalgebras of $\so(5,5)$} \\
\hline
\endhead
\hline
\endfoot
\endlastfoot

\hline
type & real subalgebras & \\
\hline

$A_4+T_1$ & $\ssl(5,\R)\oplus\so(1,1)$ & R \\
   \hline

$A_1+A_1+A_3$ & $
   \begin{aligned}
\ssl(2,\R)\oplus \ssl(2,\R)\oplus \ssl(4,\R)\\
\ssl(2,\C)\oplus \su(2,2) (2\times)\\
\su(2)\oplus \su(2)\oplus \ssl(2,\HH) (2\times)\\
   \end{aligned} $ & R \\
   \hline

$D_4+T_1$ & $
   \begin{aligned}
\so(4,4)\oplus\so(1,1)\\
\so(3,5)\oplus\u(1) (2\times)\\
   \end{aligned} $ & R \\
   \hline

$B_2+B_2$ & $
   \begin{aligned}
\so(5)\oplus \so(5)\\
\so(2,3)\oplus \so(2,3)\\
\so(4,1)\oplus \so(4,1)\\
   \end{aligned} $ & S \\
   \hline

$A_1+B_3$ & $
   \begin{aligned}
\su(2)\oplus \so(2,5)\\
\ssl(2,\R)\oplus \so(4,3)\\
   \end{aligned} $ & S \\
   \hline

$B_4$ & $\so(4,5)$ & S \\
   \hline

\end{longtable}

\newpage

\noindent{\bf Type $A_6$}

\begin{longtable}{|l|c|c|}
\caption{Maximal subalgebras of $\su(1,6)$.}\label{tab:su16}
\endfirsthead
\hline
\multicolumn{3}{|l|}{\small\slshape subalgebras of $\su(1,6)$} \\
\hline
\endhead
\hline
\endfoot
\endlastfoot

\hline
type & real subalgebras & \\
\hline

$A_2+A_3+T_1$ & $
   \begin{aligned}
\su(1,2)\oplus \su(4)\oplus\u(1)\\
\su(3)\oplus \su(1,3)\oplus\u(1)\\
   \end{aligned} $ & R \\
   \hline

$A_1+A_4+T_1$ & $
   \begin{aligned}
\su(2)\oplus \su(1,4)\oplus\u(1)\\
\ssl(2,\R)\oplus \su(5)\oplus\u(1)\\
   \end{aligned} $ & R \\
   \hline

$A_5+T_1$ & $
   \begin{aligned}
\su(1,5)\oplus\u(1)\\
\su(6)\oplus\u(1)\\
   \end{aligned} $ & R \\
   \hline

$B_3$ & $\so(6,1)$ & S \\
   \hline

\end{longtable}

\begin{longtable}{|l|c|c|}
\caption{Maximal subalgebras of $\su(2,5)$.}\label{tab:su25}
\endfirsthead
\hline
\multicolumn{3}{|l|}{\small\slshape subalgebras of $\su(2,5)$} \\
\hline
\endhead
\hline
\endfoot
\endlastfoot

\hline
type & real subalgebras & \\
\hline

$A_2+A_3+T_1$ & $
   \begin{aligned}
\su(1,2)\oplus \su(1,3)\oplus\u(1)\\
\su(3)\oplus \su(2,2)\oplus\u(1)\\
\su(1,2)\oplus \su(4)\oplus\u(1)\\
   \end{aligned} $ & R \\
   \hline

$A_1+A_4+T_1$ & $
   \begin{aligned}
\su(2)\oplus \su(2,3)\oplus\u(1)\\
\ssl(2,\R)\oplus \su(1,4)\oplus\u(1)\\
\su(2)\oplus \su(5)\oplus\u(1)\\
   \end{aligned} $ & R \\
   \hline

$A_5+T_1$ & $
   \begin{aligned}
\su(2,4)\oplus\u(1)\\
\su(1,5)\oplus\u(1)\\
   \end{aligned} $ & R \\
   \hline

$B_3$ & $\so(2,5)$ & S\\
   \hline

\end{longtable}

\begin{longtable}{|l|c|c|}
\caption{Maximal subalgebras of $\su(3,4)$.}\label{tab:su34}
\endfirsthead
\hline
\multicolumn{3}{|l|}{\small\slshape subalgebras of $\su(3,4)$} \\
\hline
\endhead
\hline
\endfoot
\endlastfoot

\hline
type & real subalgebras & \\
\hline

$A_2+A_3+T_1$ & $
   \begin{aligned}
\su(1,2)\oplus \su(2,2)\oplus\u(1)\\
\su(1,2)\oplus \su(1,3)\oplus\u(1)\\
\su(3)\oplus \su(1,3)\oplus\u(1)\\
\su(3)\oplus \su(4)\oplus\u(1)\\
   \end{aligned} $ & R \\
   \hline

$A_1+A_4+T_1$ & $
   \begin{aligned}
\ssl(2,\R)\oplus \su(2,3)\oplus\u(1)\\
\su(2)\oplus \su(2,3)\oplus\u(1)\\
\su(2)\oplus \su(1,4)\oplus\u(1)\\
   \end{aligned} $ & R \\
   \hline

$A_5+T_1$ & $
   \begin{aligned}
\su(3,3)\oplus\u(1)\\
\su(2,4)\oplus\u(1)\\
   \end{aligned} $ & R \\
   \hline

$B_3$ & $\so(4,3)$ & S \\
   \hline

\end{longtable}

\begin{longtable}{|l|c|c|}
\caption{Maximal subalgebras of $\ssl(7,\R)$.}\label{tab:sl7R}
\endfirsthead
\hline
\multicolumn{3}{|l|}{\small\slshape subalgebras of $\ssl(7,\R)$} \\
\hline
\endhead
\hline
\endfoot
\endlastfoot

\hline
type & real subalgebras & \\
\hline

$A_2+A_3+T_1$ & $\ssl(3,\R)\oplus \ssl(4,\R)\oplus\so(1,1)$ & R\\
   \hline

$A_1+A_4+T_1$ & $\ssl(2,\R)\oplus \ssl(5,\R)\oplus\so(1,1)$ & R \\
   \hline

$A_5+T_1$ & $\ssl(6,\R)\oplus\so(1,1)$ & R \\
   \hline

$B_3$ & $
   \begin{aligned}
\so(7)\\
\so(2,5)\\
\so(4,3)\\
\so(6,1)\\
   \end{aligned} $ & S \\
   \hline

\end{longtable}

\noindent{\bf Type $B_6$}

\begin{longtable}{|l|c|c|}
\caption{Maximal subalgebras of $\so(2,11)$.}\label{tab:so211}
\endfirsthead
\hline
\multicolumn{3}{|l|}{\small\slshape subalgebras of $\so(2,11)$} \\
\hline
\endhead
\hline
\endfoot
\endlastfoot

\hline
type & real subalgebras & \\
\hline

$A_3+B_3$ & $
   \begin{aligned}
\su(2,2)\oplus \so(7)\\
\ssl(2,\HH)\oplus \so(6,1)\\
\su(4)\oplus \so(2,5)\\
   \end{aligned} $ & R \\
   \hline

$B_2+D_4$ & $
   \begin{aligned}
\so(5)\oplus \so^*(8)\\
\so(4,1)\oplus \so(1,7)\\
\so(2,3)\oplus \so(8)\\
   \end{aligned} $ & R \\
   \hline

$A_1+A_1+B_4$ & $
   \begin{aligned}
\ssl(2,\R)\oplus \ssl(2,\R)\oplus \so(9)\\
\ssl(2,\C)\oplus \so(8,1)\\
\su(2)\oplus \su(2)\oplus \so(2,7)\\
   \end{aligned} $ & R \\
   \hline

$A_1+D_5$ & $
   \begin{aligned}
\su(2)\oplus \so(2,8)\\
\ssl(2,\R)\oplus \so(9,1)\\
\ssl(2,\R)\oplus \so(10)\\
   \end{aligned} $ & R \\
   \hline

$D_6$ & $
   \begin{aligned}
\so(2,10)\\
\so(11,1)\\
   \end{aligned} $ & R \\
   \hline

$B_5+T_1$ & $
   \begin{aligned}
\so(10,1)\oplus\so(1,1)\\
\so(2,9)\oplus\u(1)\\
\so(11)\oplus\u(1)\\
   \end{aligned} $ & R \\
   \hline

\end{longtable}

\begin{longtable}{|l|c|c|}
\caption{Maximal subalgebras of $\so(4,9)$.}\label{tab:so49}
\endfirsthead
\hline
\multicolumn{3}{|l|}{\small\slshape subalgebras of $\so(4,9)$} \\
\hline
\endhead
\hline
\endfoot
\endlastfoot

\hline
type & real subalgebras & \\
\hline

$A_3+B_3$ & $
   \begin{aligned}
\ssl(4,\R)\oplus \so(6,1)\\
\su(2,2)\oplus \so(2,5)\\
\ssl(2,\HH)\oplus \so(4,3)\\
\su(4)\oplus \so(4,3)\\
\su(2,2)\oplus \so(7)\\
   \end{aligned} $ & R \\
   \hline

$B_2+D_4$ & $
   \begin{aligned}
\so(5)\oplus \so(4,4)\\
\so(4,1)\oplus \so(3,5)\\
\so(2,3)\oplus \so^*(8)\\
\so(2,3)\oplus \so(1,7)\\
\so(4,1)\oplus \so(8)\\
   \end{aligned} $ & R \\
   \hline

$A_1+A_1+B_4$ & $
   \begin{aligned}
\ssl(2,\R)\oplus \ssl(2,\R)\oplus \so(2,7)\\
\ssl(2,\C)\oplus \so(6,3)\\
\su(2)\oplus \su(2)\oplus \so(4,5)\\
\ssl(2,\C)\oplus \so(8,1)\\
\su(2)\oplus \su(2)\oplus \so(9)\\
   \end{aligned} $ & R \\
   \hline

$A_1+D_5$ & $
   \begin{aligned}
\su(2)\oplus \so(4,6)\\
\ssl(2,\R)\oplus \so(3,7)\\
\ssl(2,\R)\oplus \so(2,8)\\
\su(2)\oplus \so(9,1)\\
   \end{aligned} $ & R \\
   \hline

$D_6$ & $
   \begin{aligned}
\so(4,8)\\
\so(3,9)\\
   \end{aligned} $ & R \\
   \hline

$B_5+T_1$ & $
   \begin{aligned}
\so(8,3)\oplus\so(1,1)\\
\so(4,7)\oplus\u(1)\\
\so(2,9)\oplus\u(1)\\
   \end{aligned} $ & R \\
   \hline

\end{longtable}

\begin{longtable}{|l|c|c|}
\caption{Maximal subalgebras of $\so(6,7)$.}\label{tab:so67}
\endfirsthead
\hline
\multicolumn{3}{|l|}{\small\slshape subalgebras of $\so(6,7)$} \\
\hline
\endhead
\hline
\endfoot
\endlastfoot

\hline
type & real subalgebras & \\
\hline

$A_3+B_3$ & $
   \begin{aligned}
\ssl(4,\R)\oplus \so(4,3)\\
\su(2,2)\oplus \so(4,3)\\
\su(2,2)\oplus \so(2,5)\\
\ssl(2,\HH)\oplus \so(2,5)\\
\ssl(2,\HH)\oplus \so(6,1)\\
\su(4)\oplus \so(6,1)\\
\su(4)\oplus \so(7)\\
   \end{aligned} $ & R \\
   \hline

$B_2+D_4$ & $
   \begin{aligned}
\so(2,3)\oplus \so(4,4)\\
\so(2,3)\oplus \so(3,5)\\
\so(4,1)\oplus \so(3,5)\\
\so(4,1)\oplus \so^*(8)\\
\so(5)\oplus \so^*(8)\\
\so(5)\oplus \so(1,7)\\
   \end{aligned} $ & R \\
   \hline

$A_1+A_1+B_4$ & $
   \begin{aligned}
\ssl(2,\R)\oplus \ssl(2,\R)\oplus \so(4,5)\\
\ssl(2,\C)\oplus \so(4,5)\\
\ssl(2,\C)\oplus \so(6,3)\\
\su(2)\oplus \su(2)\oplus \so(6,3)\\
\su(2)\oplus \su(2)\oplus \so(2,7)\\
   \end{aligned} $ & R \\
   \hline

$A_1+D_5$ & $
   \begin{aligned}
\ssl(2,\R)\oplus \so(5,5)\\
\ssl(2,\R)\oplus \so(4,6)\\
\su(2)\oplus \so(4,6)\\
\su(2)\oplus \so(3,7)\\
   \end{aligned} $ & R \\
   \hline

$D_6$ & $
   \begin{aligned}
\so(6,6)\\
\so(5,7)\\
   \end{aligned} $ & R \\
   \hline

$B_5+T_1$ & $
   \begin{aligned}
\so(6,5)\oplus\so(1,1)\\
\so(6,5)\oplus\u(1)\\
\so(4,7)\oplus\u(1)\\
   \end{aligned} $ & R \\
   \hline

$A_1$ & $\ssl(2,\R)$ & S \\
   \hline

\end{longtable}

\newpage

\begin{longtable}{|l|c|c|}
\caption{Maximal subalgebras of $\so(8,5)$.}\label{tab:so85}
\endfirsthead
\hline
\multicolumn{3}{|l|}{\small\slshape subalgebras of $\so(8,5)$} \\
\hline
\endhead
\hline
\endfoot
\endlastfoot

\hline
type & real subalgebras & \\
\hline

$A_3+B_3$ & $
   \begin{aligned}
\ssl(4,\R)\oplus \so(2,5)\\
\su(2,2)\oplus \so(4,3)\\
\ssl(2,\HH)\oplus \so(4,3)\\
\su(2,2)\oplus \so(6,1)\\
\su(4)\oplus \so(2,5)\\
\ssl(2,\HH)\oplus \so(7)\\
   \end{aligned} $ & R \\
   \hline

$B_2+D_4$ & $
   \begin{aligned}
\so(4,1)\oplus \so(4,4)\\
\so(2,3)\oplus \so(3,5)\\
\so(2,3)\oplus \so^*(8)\\
\so(5)\oplus \so(3,5)\\
\so(4,1)\oplus \so(1,7)\\
\so(5)\oplus \so(8)\\
   \end{aligned} $ & R \\
   \hline

$A_1+A_1+B_4$ & $
   \begin{aligned}
\ssl(2,\R)\oplus \ssl(2,\R)\oplus \so(6,3)\\
\ssl(2,\C)\oplus \so(4,5)\\
\su(2)\oplus \su(2)\oplus \so(4,5)\\
\ssl(2,\C)\oplus \so(2,7)\\
\su(2)\oplus \su(2)\oplus \so(8,1)\\
   \end{aligned} $ & R \\
   \hline

$A_1+D_5$ & $
   \begin{aligned}
\su(2)\oplus \so(5,5)\\
\ssl(2,\R)\oplus \so(4,6)\\
\ssl(2,\R)\oplus \so(3,7)\\
\su(2)\oplus \so(2,8)\\
   \end{aligned} $ & R \\
   \hline

$D_6$ & $
   \begin{aligned}
\so(5,7)\\
\so(4,8)\\
   \end{aligned} $ & R \\
   \hline

$B_5+T_1$ & $
   \begin{aligned}
\so(4,7)\oplus\so(1,1)\\
\so(6,5)\oplus\u(1)\\
\so(8,3)\oplus\u(1)\\
   \end{aligned} $ & R \\
   \hline

\end{longtable}

\begin{longtable}{|l|c|c|}
\caption{Maximal subalgebras of $\so(10,3)$.}\label{tab:so103}
\endfirsthead
\hline
\multicolumn{3}{|l|}{\small\slshape subalgebras of $\so(10,3)$} \\
\hline
\endhead
\hline
\endfoot
\endlastfoot

\hline
type & real subalgebras & \\
\hline

$A_3+B_3$ & $
   \begin{aligned}
\ssl(4,\R)\oplus \so(7)\\
\su(2,2)\oplus \so(6,1)\\
\ssl(2,\HH)\oplus \so(2,5)\\
\su(4)\oplus \so(4,3)\\
   \end{aligned} $ & R \\
   \hline

$B_2+D_4$ & $
   \begin{aligned}
\so(5)\oplus \so(3,5)\\
\so(4,1)\oplus \so^*(8)\\
\so(2,3)\oplus \so(1,7)\\
\so(2,3)\oplus \so(8)\\
   \end{aligned} $ & R \\
   \hline

$A_1+A_1+B_4$ & $
   \begin{aligned}
\ssl(2,\R)\oplus \ssl(2,\R)\oplus \so(8,1)\\
\ssl(2,\C)\oplus \so(2,7)\\
\su(2)\oplus \su(2)\oplus \so(6,3)\\
\ssl(2,\C)\oplus \so(9)\\
   \end{aligned} $ & R \\
   \hline

$A_1+D_5$ & $
   \begin{aligned}
\su(2)\oplus \so(3,7)\\
\ssl(2,\R)\oplus \so(2,8)\\
\ssl(2,\R)\oplus \so(9,1)\\
\su(2)\oplus \so(10)\\
   \end{aligned} $ & R \\
   \hline

$D_6$ & $
   \begin{aligned}
\so(3,9)\\
\so(2,10)\\
   \end{aligned} $ & R \\
   \hline

$B_5+T_1$ & $
   \begin{aligned}
\so(2,9)\oplus\so(1,1)\\
\so(8,3)\oplus\u(1)\\
\so(10,1)\oplus\u(1)\\
   \end{aligned} $ & R \\
   \hline

\end{longtable}

\begin{longtable}{|l|c|c|}
\caption{Maximal subalgebras of $\so(12,1)$.}\label{tab:so121}
\endfirsthead
\hline
\multicolumn{3}{|l|}{\small\slshape subalgebras of $\so(12,1)$} \\
\hline
\endhead
\hline
\endfoot
\endlastfoot

\hline
type & real subalgebras & \\
\hline

$A_3+B_3$ & $
   \begin{aligned}
\ssl(2,\HH)\oplus \so(7)\\
\su(4)\oplus \so(6,1)\\
   \end{aligned} $ & R \\
   \hline

$B_2+D_4$ & $
   \begin{aligned}
\so(5)\oplus \so(1,7)\\
\so(4,1)\oplus \so(8)\\
   \end{aligned} $ & R \\
   \hline

$A_1+A_1+B_4$ & $
   \begin{aligned}
\ssl(2,\C)\oplus \so(9)\\
\su(2)\oplus \su(2)\oplus \so(8,1)\\
   \end{aligned} $ & R \\
   \hline

$A_1+D_5$ & $
   \begin{aligned}
\su(2)\oplus \so(9,1)\\
\ssl(2,\R)\oplus \so(10)\\
   \end{aligned} $ & R \\
   \hline

$D_6$ & $
   \begin{aligned}
\so(11,1)\\
\so(12)\\
   \end{aligned} $ & R \\
   \hline

$B_5+T_1$ & $
   \begin{aligned}
\so(11)\oplus\so(1,1)\\
\so(10,1)\oplus\u(1)\\
   \end{aligned} $ & R \\
   \hline

\end{longtable}

\noindent{\bf Type $C_6$}

\begin{longtable}{|l|c|c|}
\caption{Maximal subalgebras of $\sp(1,5)$.}\label{tab:sp15}
\endfirsthead
\hline
\multicolumn{3}{|l|}{\small\slshape subalgebras of $\sp(1,5)$} \\
\hline
\endhead
\hline
\endfoot
\endlastfoot

\hline
type & real subalgebras & \\
\hline

$A_5+T_1$ & $\su(1,5)\oplus\u(1)$ & R \\
   \hline

$C_3+C_3$ & $\sp(3)\oplus \sp(1,2)$ & R \\
   \hline

$B_2+C_4$ & $
   \begin{aligned}
\so(4,1)\oplus \sp(4)\\
\so(5)\oplus \sp(1,3)\\
   \end{aligned} $ & R \\
   \hline

$A_1+C_5$ & $
   \begin{aligned}
\su(2)\oplus \sp(1,4)\\
\su(2)\oplus \sp(5)\\
   \end{aligned} $ & R \\
   \hline

$A_1+A_3$ & $\su(2)\oplus \ssl(2,\HH)$ & S \\
   \hline

\end{longtable}

\begin{longtable}{|l|c|c|}
\caption{Maximal subalgebras of $\sp(2,4)$.}\label{tab:sp24}
\endfirsthead
\hline
\multicolumn{3}{|l|}{\small\slshape subalgebras of $\sp(2,4)$} \\
\hline
\endhead
\hline
\endfoot
\endlastfoot

\hline
type & real subalgebras & \\
\hline

$A_5+T_1$ & $\su(2,4)\oplus\u(1)$ & R \\
   \hline

$C_3+C_3$ & $
   \begin{aligned}
\sp(1,2)\oplus \sp(1,2)\\
\sp(3)\oplus \sp(1,2)\\
   \end{aligned} $ & R \\
   \hline

$B_2+C_4$ & $
   \begin{aligned}
\so(4,1)\oplus \sp(1,3)\\
\so(5)\oplus \sp(2,2)\\
\so(5)\oplus \sp(4)\\
   \end{aligned} $ & R \\
   \hline

$A_1+C_5$ & $
   \begin{aligned}
\su(2)\oplus \sp(2,3)\\
\su(2)\oplus \sp(1,4)\\
   \end{aligned} $ & R \\
   \hline

$A_1+B_2$ & $\ssl(2,\R)\oplus \so(5)$ & S\\
   \hline

$A_1+A_3$ & $\su(2)\oplus \su(2,2)$ & S \\
   \hline

\end{longtable}

\begin{longtable}{|l|c|c|}
\caption{Maximal subalgebras of $\sp(3,3)$.}\label{tab:sp33}
\endfirsthead
\hline
\multicolumn{3}{|l|}{\small\slshape subalgebras of $\sp(3,3)$} \\
\hline
\endhead
\hline
\endfoot
\endlastfoot

\hline
type & real subalgebras & \\
\hline

$A_5+T_1$ & $
   \begin{aligned}
\ssl(3,\HH)\oplus\so(1,1)\\
\su(3,3)\oplus\u(1)\\
   \end{aligned} $ & R \\
   \hline

$C_3+C_3$ & $
   \begin{aligned}
\sp(3,\C)\\
\sp(1,2)\oplus \sp(1,2)\\
\sp(3)\oplus \sp(3)\\
   \end{aligned} $ & R \\
   \hline

$B_2+C_4$ & $
   \begin{aligned}
\so(4,1)\oplus \sp(2,2)\\
\so(5)\oplus \sp(1,3) (2\times)\\
   \end{aligned} $ & R \\
   \hline

$A_1+C_5$ & $\su(2)\oplus \sp(2,3) (2\times)$ & R\\
   \hline

$A_1+B_2$ & $\su(2)+\so(4,1)$ & S \\
   \hline

$A_1+A_3$ & $
   \begin{aligned}
     \su(2)\oplus \ssl(4,\R)\\
     \ssl(2,\R)\oplus \su(1,3)\\
   \end{aligned} $ & S \\
   \hline

\end{longtable}%

\newcommand{\etalchar}[1]{$^{#1}$}

\end{document}